\documentclass[titlepage,12pt]{article} 
\usepackage{hyperref}
\usepackage[usenames,dvipsnames]{xcolor}
\usepackage[title]{appendix}
\usepackage{amssymb,amsthm,amsmath} 
\usepackage[a4paper]{geometry}
\usepackage{datetime2}
\usepackage{enumitem}
\usepackage[utf8]{inputenc}

\date{}

\newcommand{\boxa}[1]{\fbox{$\displaystyle#1$}}

\newcommand{\dath}{D(B)\times H}
\newcommand{\hg}{\widehat{g}}
\newcommand{\gaph}{\mathcal{G}(H,\lambda_{1},\lambda_{2})}

\newcommand{\ep}{\varepsilon}
\renewcommand{\qed}{{\penalty 10000\mbox{$\quad\Box$}}}
\newcommand{\re}{\mathbb{R}}

\newtheorem{thm}{Theorem}[section]

\newtheorem{rmk}[thm]{Remark}
\newtheorem{prop}[thm]{Proposition}
\newtheorem{defn}[thm]{Definition}
\newtheorem{cor}[thm]{Corollary}

\newtheorem{lemma}[thm]{Lemma}

\title{Small perturbations for a Duffing-like evolution equation involving non-commuting operators}

\author{Marina Ghisi\vspace{1ex}\\ 
{\normalsize Universit\`a degli Studi di Pisa} \\
{\normalsize Dipartimento di Matematica}\\ 
{\normalsize PISA (Italy)}\\
{\normalsize e-mail: \texttt{marina.ghisi@unipi.it}}
\and
Massimo Gobbino\vspace{1ex}\\ 
{\normalsize Universit\`a degli Studi di Pisa} \\
{\normalsize Dipartimento di Ingegneria Civile e Industriale}\\ 
{\normalsize PISA (Italy)}\\  
{\normalsize e-mail: \texttt{massimo.gobbino@unipi.it}}
\and
Alain Haraux\vspace{1ex}\\ 
{\normalsize Sorbonne Universit\'e, Universit\'e Paris-Diderot SPC, CNRS, INRIA}, \\
{\normalsize Laboratoire Jacques-Louis Lions,  LJLL, F-75005,
Paris, France.}\\ 
{\normalsize e-mail: \texttt{haraux@ann.jussieu.fr}}}

\begin{document}
\maketitle

\begin{abstract}

We consider an abstract evolution equation with linear damping, a nonlinear term of Duffing type, and a small forcing term. The abstract problem is inspired by some models for damped oscillations of a beam subject to external loads or magnetic fields, and shaken by a transversal force. 

The main feature is that very natural choices of the boundary conditions lead to equations whose linear part involves two operators that do not commute.

We extend to this setting the results that are known in the commutative case, namely that for asymptotically small forcing terms all solutions are eventually close to the three equilibrium points of the unforced equation, two stable and one unstable.
\vspace{6ex}

\noindent{\bf Mathematics Subject Classification 2010 (MSC2010):} 
35B40, 35L75, 35L90.

		
\vspace{6ex}


\end{abstract}

 
\section{Introduction}

Let us consider the partial differential equation
\begin{equation}
u_{tt}+\delta u_{t}+k_{1}u_{xxxx}+k_{2}u_{xx}-k_{3}\left(\int_{0}^{1}u_{x}^{2}\,dx\right)u_{xx}=f(t, x)
\label{eqn:duffing-concr}
\end{equation}
in the strip $(t,x)\in [0,+\infty)\times [0,1]$, where $\delta$, $k_{1}$, $k_{2}$, $k_{3}$ are positive constants, and $f:[0,+\infty)\times[0,1]\to\re$ is a given function (forcing term).

Equation (\ref{eqn:duffing-concr}) was derived as a model for the motion of a beam in different physical systems, for example
\begin{itemize}

\item  in~\cite{holm-mars} the beam is buckled by an external load $k_{2}$, and shaken by a transverse displacement (depending only on time, in that model),

\item  in~\cite{moon-holm} (the so-called magneto-elastic cantilever beam) the beam is clamped vertically at the upper end, and suspended at the other end between two magnets secured to a base, and the whole system is shaken by an external force transversal to the beam. 

\end{itemize}

Equation (\ref{eqn:duffing-concr}) may be seen as an abstract evolution problem in a Hilbert space, but the precise setting depends on the boundary conditions. 

\paragraph{\textmd{\textit{The ``commutative'' case}}} 

Let us consider equation (\ref{eqn:duffing-concr}) with boundary conditions
\begin{equation}
u(t,x)=u_{xx}(t,x)=0
\qquad
\forall(t,x)\in [0,+\infty)\times \{0,1\},
\label{eqn:bc-hinged}
\end{equation}
physically corresponding to ``hinged ends''. In this case (\ref{eqn:duffing-concr}) may be seen as an abstract evolution problem of the form
\begin{equation}
u''+\delta u'+k_{1}A^{2}u-k_{2}Au+k_{3}|A^{1/2}u|^{2}Au=f(t)
\label{eqn:duffing-AAk}
\end{equation}
in the Hilbert space $H:=L^{2}((0,1))$, where $Au=-u_{xx}$ with domain 
\begin{equation}
D(A):=\left\{u\in H^{2}((0,1)): u(0)=u(1)=0\right\}.
\nonumber
\end{equation} 

Up to changing the unknown and the operator according to the rules
$$u(t)\rightsquigarrow\alpha u(\beta t),
\qquad\qquad
A\rightsquigarrow\gamma A$$
for suitable values of $\alpha$, $\beta$, $\gamma$, we can assume that three of the four constants in (\ref{eqn:duffing-AAk}) are equal to~1, and we end up with an equation of the form 
\begin{equation}
u''+u'+A^{2}u-\lambda Au+|A^{1/2}u|^{2}Au=f(t),
\label{eqn:duffing-AA}
\end{equation}
depending only on one positive parameter $\lambda$. We point out that with these choices it turns out that, up to constants, $A^{2}u=u_{xxxx}$ with domain
\begin{equation}
D(A^{2}):=\left\{u\in H^{4}((0,1)):u(0)=u_{xx}(0)=u(1)=u_{xx}(1)=0\right\}.
\nonumber
\end{equation}

Equation \eqref {eqn:duffing-AA}  can be considered more generally whenever   $A$ is  a coercive  selfadjoint operator and $ A^{-1}:  H \rightarrow H $  is compact. In this case $H$ admits a countable orthonormal system $\{e_{n}\}$ made of  the eigenvectors of $A$. The theory  has been done when the first eigenvalues $\lambda_{1}$  of $A$ is simple. The behavior of solutions to (\ref{eqn:duffing-AA}) depends on the position of $\lambda$ with respect to the eigenvalues of $A$. When $\lambda<\lambda_{1}$ the operator $A^{2}-\lambda A$ is positive, the functional
\begin{equation}
\mathcal{E}_{A,A}(u):=\frac{1}{2}|Au|^{2}-\frac{\lambda}{2}|A^{1/2}u|^{2}+\frac{1}{4}|A^{1/2}u|^{4}
\label{defn:energy-landscape}
\end{equation}
is convex and has a unique minimum point at the origin, and the trivial solution $u(t)\equiv 0$ is the unique stationary solution of equation (\ref{eqn:duffing-AA}) in the unforced case $f(t)\equiv 0$. If $f(t)$ is asymptotically small enough, then all solutions are asymptotic to each other as $t\to +\infty$, and lie eventually in a neighborhood of the origin whose radius depends on the asymptotic size of the forcing term. We refer to~\cite{Aloui-Haraux,F-H-DCDS,F-H-JMPA,loud2} for significant results in the convex case 

The case with $\lambda\in(\lambda_{1},\lambda_{2})$ was investigated in~\cite{GGH:duffing}. Now the operator $A^{2}-\lambda A$ is negative in the direction spanned by $e_{1}$, and positive in the orthogonal space. The functional $\mathcal{E}_{A,A}(u)$ has three stationary points: the origin, which is no longer a minimum point, and two minimum points of the form $\pm\sigma_{0}e_{1}$, with $\sigma_{0}=(\lambda-\lambda_{1})^{1/2}\lambda_{1}^{-1/2}$. 

As a consequence, in the unforced case $f(t)\equiv 0$ equation (\ref{eqn:duffing-concr}) has three stationary solutions: the trivial solution $u(t)\equiv 0$, which is now unstable, and the two stable solutions of the form $u(t)\equiv\pm\sigma_{0}e_{1}$, corresponding to the minimum points of the functional $\mathcal{E}_{A,A}(u)$. In the forced case with an external force that is asymptotically small enough, all solutions fall eventually in a neighborhood of one of the three stationary points, within a distance depending on the asymptotic size of the forcing, and any two solutions that are eventually close to the same stationary point are actually asymptotic to each other.

When $\lambda>\lambda_{2}$, the number of stationary points of the functional increases, as well as the number of stationary solutions to (\ref{eqn:duffing-AA}) in the unforced case. This regime has not been investigated explicitly, but the same approach as in~\cite{GGH:duffing} is likely to work when all eigenvalues are simple, i.e. $0<\lambda_{1}<\lambda_{2}<\ldots$ or more generally as long as we do not cross a multiple eigenvalue.

We conclude this paragraph by mentioning two more sets of boundary conditions that lead to commutative operators.
\begin{itemize}

\item  The periodic boundary conditions
$$u(t,0)=u(t,1),
\qquad
u_{x}(t,0)=u_{x}(t,1),
\qquad
\forall t\geq 0,$$
$$u_{xx}(t,0)=u_{xx}(t,1),
\qquad
u_{xxx}(t,0)=u_{xxx}(t,1),
\qquad
\forall t\geq 0,$$
in which case the operator $A$ acts again as $Au=-u_{xx}$, but now with domain
\begin{equation}
D(A):=\left\{u\in H^{2}((0,1)):u(0)=u(1),\ u_{x}(0)=u_{x}(1)\right\}.
\nonumber
\end{equation}

\item  Boundary conditions such as
$$u(t,0)=u_{xx}(t,0)=0,
\qquad
u_{x}(t,1)=u_{xxx}(t,1)=0
\qquad
\forall t\geq 0.$$

Indeed such a case can be easily reduced to (\ref{eqn:bc-hinged}) after extending the solution to the interval $(0,2)$ by means of a reflection with respect to $x=1$.

\end{itemize}

\paragraph{\textmd{\textit{The ``non-commutative'' case}}} 

Let us consider now equation (\ref{eqn:duffing-concr}) with boundary conditions
\begin{equation}
u(t,x)=u_{x}(t,x)=0
\qquad
\forall(t,x)\in [0,+\infty)\times \{0,1\},
\nonumber
\end{equation}
physically corresponding to ``clamped ends''. After suitable variable changes, we end up with an abstract evolution problem of the form
\begin{equation}
u''+u'+B^{2}u-\lambda Au+|A^{1/2}u|^{2}Au=f(t).
\label{eqn:duffing-AB}
\end{equation}

The Hilbert space $H$ and the operator $A$ are the same as before. Also the operator $B^{2}$ acts as $B^{2}u=u_{xxxx}$ as in the previous case, but now with domain
\begin{equation}
D(B^{2}):=\left\{u\in H^{4}((0,1)):u(0)=u_{x}(0)=u(1)=u_{x}(1)=0\right\}.
\nonumber
\end{equation}

This makes a great difference, because $A^{2}$ and $B^{2}$ have now different eigenspaces, and hence they do not commute (note also that, with this choice of the domain $D(B^{2})$, the operator $B$, defined as the square root of $B^{2}$, does not act as $-u_{xx}$).

Nevertheless, the functional has now the form
\begin{equation}
\mathcal{E}_{A,B}(u):=\frac{1}{2}|Bu|^{2}-\frac{\lambda}{2}|A^{1/2}u|^{2}+\frac{1}{4}|A^{1/2}u|^{4},
\label{defn:funct-AB}
\end{equation}
which is qualitatively similar to (\ref{defn:energy-landscape}). In particular, in Proposition~\ref{prop:gaph} we show that there exist again two positive constants $\lambda_{2}>\lambda_{1}>0$, which are now the two smallest eigenvalues of the operator $A^{-1}B^{2}$ (see Proposition~\ref{prop:concrete} and the final appendix), with the following properties.

\begin{itemize}

\item  When $\lambda<\lambda_{1}$ the operator $B^{2}-\lambda A$ is positive, and the functional $\mathcal{E}_{A,B}(u)$ is convex with a unique minimum point at the origin.

\item  When $\lambda_{1}<\lambda<\lambda_{2}$ the operator $B^{2}-\lambda A$ is negative in a subspace of dimension one, and positive in the orthogonal subspace. In this regime the functional $\mathcal{E}_{A,B}(u)$ has three stationary points: the origin, which in no longer a minimum point, and two minimum points that are symmetric with respect to the origin.

\item  When $\lambda>\lambda_{2}$ the operator $B^{2}-\lambda A$ is negative in a subspace of dimension at least two, and the functional $\mathcal{E}_{A,B}(u)$ has more than three stationary points.

\end{itemize}

In this paper we investigate the regime $\lambda\in(\lambda_{1},\lambda_{2})$, and we show that solutions to (\ref{eqn:duffing-AB}) have the same qualitative behavior as the solutions to the ``commutative'' model (\ref{eqn:duffing-AA}) in the corresponding regime.

We conclude this paragraph by mentioning that the non-commutative case is also the correct setting for dealing with boundary conditions such as
$$u(t,0)=u_{xx}(t,0)=u(t,1)=u_{x}(t,1)=0
\qquad
\forall t\geq 0,$$
physically corresponding to a beam hinged in $x=0$, and clamped in $x=1$. In this case $A$ is the same operator as before, and $B^{2}u=u_{xxxx}$ with domain
$$D(B^{2}):=\left\{u\in H^{4}((0,1)):u(0)=u_{xx}(0)=u(1)=u_{x}(1)=0\right\}.$$

Unfortunately, the cantilever beam with one free end described in~\cite{moon-holm} fits neither in the commutative, nor in the non-commutative setting. That model involves  nonlinear boundary conditions in the free endpoint, and for this reason it deserves a distinct theory that we plan to investigate in the future.

\paragraph{\textmd{\textit{Structure of the paper}}} 

This paper is organized as follows. In section~\ref{sec:statements} we clarify the functional setting, we state a preliminary well-posedness result for (\ref{eqn:duffing-AB}) (Proposition~\ref{wp}), and then we state our main result (Theorem~\ref{thm:main}) concerning the existence of three different asymptotic regimes, and a simple consequence (Corollary~\ref{cor}). In section~\ref{sec:proof-sketch} we state four auxiliary propositions, where we concentrate the technical machinery of the paper. In section~\ref{sec:linear-algebra} we prove all the abstract properties of the operator $B^{2}-\lambda A$ that we need in the paper. Section~\ref{sec:proof-prop} is devoted to the proof of the four auxiliary propositions. Section~\ref{sec:proof-main} contains the proof of our main result. In section~\ref{sec:concrete} we show that the beam equation (\ref{eqn:duffing-concr}) with clamped ends fits in our abstract framework. Finally, in the appendix we discuss the correct functional setting for the operator $A^{-1}B^{2}$ in the case where $A$ and $B$ do not commute necessarily.


\setcounter{equation}{0}
\section{Statements}\label{sec:statements}

Throughout this paper we always consider equation (\ref{eqn:duffing-AB}) with initial data
\begin{equation}
u(0)=u_{0}\in D(B),
\qquad\qquad
u'(0)=u_{1}\in H.
\label{eqn:data}
\end{equation}

\paragraph{\textmd{\textit{Well-posedness}}} 

Rather classical techniques lead to the following well-posedness result under quite general assumptions on the operators $A$ and $B$, and on the parameter $\lambda$.

\begin{prop}\label{wp}

Let $H$ be a Hilbert space, let $\lambda$ be a real number, let $f:\re\to H$ be a continuous function, and let $A$ and $B$ be two self-adjoint nonnegative linear operators on $H$ with dense domains $D(B)\subseteq D(A)$.

Let us assume that there exists a positive constant $\mu_{1}$ such that 
\begin{equation}
|Bu|^{2}\geq\mu_{1}|Au|^{2}  
\qquad
\forall u\in D(B).
\label{defn:mu-1} 
\end {equation}

Then the following statements hold true.
\begin{enumerate}
\renewcommand{\labelenumi}{(\arabic{enumi})}
  
\item \emph{(Global existence and uniqueness)} For every $(u_{0},u_{1})\in\dath$, problem (\ref{eqn:duffing-AB})--(\ref{eqn:data}) admits a unique global solution
\begin{equation}
u\in C^{0}\left(\re,D(B)\right)\cap C^{1}\left(\re,H\right).
\nonumber
\end{equation}
  
\item \emph{(Continuous dependence on initial data)} Let $\{(u_{0n},u_{1n})\}\subseteq D(B)\times H$ be any sequence with
\begin{equation}
(u_{0n},u_{1n})\to(u_{0},u_{1})
\quad\quad
\mbox{in }\dath,
\nonumber
\end{equation}
and let $u_{n}(t)$ denote the solution to (\ref{eqn:duffing-AB}) with data $u_{n}(0)=u_{0n}$ and $u'(0)=u_{1n}$.
  
Then for every $T>0$ it turns out that
\begin{equation}
u_{n}(t)\to u(t)
\quad\quad
\mbox{uniformly in }C^{0}\left([-T,T],D(B)\right),
\nonumber
\end{equation}
\vspace{-\baselineskip}
\begin{equation}
u_{n}'(t)\to u'(t)
\quad\quad
\mbox{uniformly in }C^{0}\left([-T,T],H\right).
\nonumber
\end{equation}
   
\item \emph{(Derivative of the energy)} The classical energy
\begin{equation}
E(t):=\frac{1}{2}|u'(t)|^{2}+\frac{1}{2}|Bu(t)|^{2}-\frac{\lambda}{2}|A^{1/2}u(t)|^{2}+\frac{1}{4}|A^{1/2}u(t)|^{4}
\label{defn:E}
\end{equation}
is of class $C^{1}$, and its time-derivative is 
\begin{equation}
E'(t)=-|u'(t)|^{2}+\langle u'(t),f(t)\rangle
\quad\quad
\forall t\in\re.
\nonumber
\end{equation} 
\end{enumerate}

\end{prop}

\begin{rmk}
\begin{em}
For the sake of simplicity, in the statement of Proposition~\ref{wp} above we assumed that the forcing term $f(t)$ is defined for every $t\in\re$. Of course, if $f(t)$ is defined only in the half-line $t\geq 0$, or in some interval $[0,T]$, then the solution $u(t)$ is defined for the same values of $t$.
\end{em}
\end{rmk}


In the case where $A=B$, it was proved that the asymptotic dynamics depend on the position of $\lambda$ with respect to the spectrum of $A$. When the two operators $A$ and $B$ are different, and do not commute, it is not immediately clear which set will play the role of the spectrum of $A$. From the heuristic point of view, it is useful to consider first the finite dimensional case.

\paragraph{\textmd{\textit{The finite dimensional case}}} 

If $\dim H < \infty $, then the operators $A$ and $B$ are represented  by two symmetric and positive matrices, and the square roots of the quadratic forms associated to $B^{2}$ and $A$ define two equivalent norms on $H$. Moreover let us define $$ \Sigma = \{u\in H,\quad |A^{1/2}u|^{2}=1\} $$  and 
\begin{equation}
\lambda_{1}:=\min_{u\in \Sigma} |Bu|^{2}
\label{defn:lambda-1}
\end{equation}
then any minimizer $u$ of (\ref{defn:lambda-1}) satisfies
\begin{equation}
B^{2}u=\lambda_{1}Au.
\end{equation}

In particular, $\lambda_{1}$ is the smallest eigenvalue of the matrix $A^{-1}B^{2}$, and the set of minimizers of (\ref{defn:lambda-1}) spans the corresponding eigenspace. From the definition of $\lambda_{1}$ it follows also that the matrix $B^{2}-\lambda A$ is positive for every $\lambda<\lambda_{1}$.

Now let us choose a minimizer $e_{1}\in \Sigma$ of (\ref{defn:lambda-1}) and let us set $$ \Sigma_+ = \{u\in H,\quad |A^{1/2}u|^{2}=1,\ \langle u,Ae_{1}\rangle=0\} $$  and 
\begin{equation}
\lambda_{2}:=\min_{u\in \Sigma_+} |Bu|^{2}\geq\lambda_{1}.
\label{defn:lambda-2}
\end{equation}

Then it turns out that $\lambda_{2}$ is the second smallest eigenvalue of $A^{-1}B^{2}$, and the set of minimizers of (\ref{defn:lambda-2}) is the corresponding eigenspace. If the strict inequality $\lambda_{2}>\lambda_{1}$ holds true (and this happens if and only if $\lambda_{1}$ is simple), then for every $\lambda\in(\lambda_{1},\lambda_{2})$ the matrix $B^{2}-\lambda A$ has exactly one negative eigenvalue, while all remaining eigenvalues are positive. In this case one says that the negative inertia index of $B^{2}-\lambda A$ is 1. 

For $\lambda>\lambda_{2}$, the matrix $B^{2}-\lambda A$ has at least two negative eigenvalues.

This process can be carried on, thus showing that the number of negative eigenvalues of $B^{2}-\lambda A$ increases when $\lambda$ crosses the eigenvalues of $A^{-1}B^{2}$.

\paragraph{\textmd{\textit{Operators with gap condition}}} 

In the following definition we extend to infinite dimensions the framework described above.

\begin{defn}[Pairs of operators with gap condition]\label{defn:gap}
\begin{em}

Let $H$ be a Hilbert space, and let $\lambda_{1}<\lambda_{2}$ be two positive real numbers. We say that two operators $A$ and $B$ satisfy the $(\lambda_{1},\lambda_{2})$ gap condition, and we write $(A,B)\in\gaph$, if 
\begin{itemize}

\item  $A$ and $B$ are self-adjoint linear operators on $H$ with dense domains $D(B)\subseteq D(A)$, and there exists a positive constant $\mu_{1}$ for which (\ref{defn:mu-1}) holds true.

\item there exists a positive real number $\mu_{2}$ such that
\begin{equation}
|A^{1/2}u|^{2}\geq\mu_{2}|u|^{2}
\qquad
\forall u\in D(A^{1/2}),
\label{defn:mu-2}
\end{equation}

\item  there exists $e_{1}\in D(B^{2})$, with $|A^{1/2} e_{1}|=1$, such that
\begin{equation}
B^{2}e_{1}=\lambda_{1}Ae_{1}
\label{hp:e1}
\end{equation}
and
\begin{equation}
|Bu|^{2}\geq\lambda_{2}|A^{1/2}u|^{2}
\qquad
\forall u\in\{x\in D(B):\langle x,Ae_{1}\rangle=0\}.
\label{hp:e1-perp}
\end{equation}

\end{itemize}

\end{em}
\end{defn}

In the following result we collect all the properties of this class of operators that we need in this paper.

\begin{prop}[Properties of pairs of operators with gap condition]\label{prop:gaph}

Let $H$ be a Hilbert space, let $\lambda_{1}<\lambda_{2}$ be two positive real numbers, and let $(A,B)\in\gaph$.

Then the following statements hold true.
\begin{enumerate}
\renewcommand{\labelenumi}{(\arabic{enumi})}

\item For every $\lambda\in \re$ the operator $B^{2}-\lambda A$ is self-adjoint as an unbounded linear operator in $H$ with domain $D(B^2)$.

\item For every $\lambda<\lambda_{1}$ the self-adjoint operator $B^{2}-\lambda A$ is positive, namely
\begin{equation}
|Bu|^{2}-\lambda|A^{1/2}u|^{2}>0
\qquad
\forall u\in D(B)\setminus\{0\}.
\label{th:prop-sub}
\end{equation}

\item  For every $\lambda\in(\lambda_{1},\lambda_{2})$ the operator $B^{2}-\lambda A$ has negative inertia index equal to one. More precisely, there exist a positive real number $\lambda_{0}$, and an element $e_{0}\in D(B^{2})$ with $|e_{0}|=1$ (both $\lambda_{0}$ and $e_{0}$ do depend also on $\lambda$), such that
\begin{equation}
B^{2}e_{0}-\lambda Ae_{0}=-\lambda_{0}e_{0},
\label{defn:e0}
\end{equation}
and there exists a positive constant $\mu_{3}$ (again depending on $\lambda$) such that
\begin{equation}
|Bu|^{2}-\lambda|A^{1/2}u|^{2}\geq\mu_{3}|Bu|^{2}
\qquad
\forall u\in\{x\in D(B):\langle x,e_{0}\rangle=0\}.
\label{defn:e0+}
\end{equation}

\item  For every $\lambda\in(\lambda_{1},\lambda_{2})$ the functional $\mathcal{E}_{A,B}(u)$ defined in (\ref{defn:funct-AB}) has three stationary points, namely the three solutions to the equation
\begin{equation}
B^{2}u-\lambda Au+|A^{1/2}u|^{2}Au=0.
\nonumber
\end{equation}

These three solutions are 0 and $\pm\sigma_{0}e_{1}$, where $e_{1}$ is the vector that appears in (\ref{hp:e1}) and (\ref{hp:e1-perp}), and
\begin{equation}
\sigma_{0}:=\sqrt{\lambda-\lambda_{1}}.
\label{defn:sigma-0}
\end{equation}

\item  If $\lambda_{2}$ is the largest constant for which (\ref{hp:e1-perp}) holds true, then for every $\lambda>\lambda_{2}$ the operator $B^{2}-\lambda A$ has negative inertia index at least two, namely there exists a two-dimensional subspace $V\subseteq D(B)$ such that
\begin{equation}
|Bv|^{2}-\lambda|A^{1/2}v|^{2}<0
\qquad
\forall v\in V\setminus\{0\}.
\label{th:prop-4}
\end{equation}

\end{enumerate}

\end{prop}


\paragraph{\textmd{\textit{Main result}}} 

From now on, we always consider problem (\ref{eqn:duffing-AB})--(\ref{eqn:data}) under the following assumptions, which we briefly call \emph{standard assumptions}:
\begin{itemize}

\item $H$ is a Hilbert space,

\item $0<\lambda_{1}<\lambda_{2}$ are real numbers,

\item $(A,B)\in\gaph$ is a pair of operators satisfying the gap condition,

\item $\lambda\in(\lambda_{1},\lambda_{2})$ is a real number,

\item  $f:[0,+\infty)\to H$ is a bounded and continuous function.

\end{itemize}

Our main result is the following.

\begin{thm}[Asymptotic behavior of solutions with small external force]\label{thm:main}

Let us consider problem (\ref{eqn:duffing-AB})--(\ref{eqn:data}) under the standard assumptions presented above. Let $\sigma_{0}$ be the constant defined by (\ref{defn:sigma-0}).

Then there exists two positive constants $\ep_{0}$ and $M_{0}$, independent of the forcing term $f(t)$ and of the solution $u(t)$, for which the following statements hold true whenever
\begin{equation}
\limsup_{t\to+\infty}|f(t)|\leq\ep_{0}.
\label{hp:main}
\end{equation}

\begin{enumerate}
\renewcommand{\labelenumi}{(\arabic{enumi})}
  
\item \emph{(Alternative)} For every solution $u(t)$ to (\ref{eqn:duffing-AB}), there exists $\sigma\in\{-\sigma_{0},0,\sigma_{0}\}$ such that
\begin{equation}
\limsup_{t\to+\infty}\left(|u'(t)|+|B(u(t)-\sigma e_{1})|\strut\right)\leq M_{0}\limsup_{t\to+\infty}|f(t)|.
\label{th:main-alternative}
\end{equation}
	   
\item \emph{(Asymptotic convergence)} If $u(t)$ and $v(t)$ are any two solutions to (\ref{eqn:duffing-AB}) satisfying (\ref{th:main-alternative}) with the same $\sigma\in\{-\sigma_{0},0,\sigma_{0}\}$, then $u(t)$ and $v(t)$ are asymptotic to each other in the sense that
\begin{equation}
\lim_{t\to +\infty}\left(|u'(t)-v'(t)|+|B(u(t)-v(t))|\strut\right)=0.
\label{th:main-asymptotic}
\end{equation}

\item \emph{(Solutions with $\sigma=\pm\sigma_{0}$)} The set of initial data $(u_{0},u_{1})$ for which the solution to (\ref{eqn:duffing-AB})--(\ref{eqn:data}) satisfies (\ref{th:main-alternative}) with a given $\sigma\in\{-\sigma_{0},\sigma_{0}\}$ is a nonempty open subset of $\dath$. 
	
\item \emph{(Solutions with $\sigma = 0$)} The set of initial data $(u_{0},u_{1})$ for which the solution to (\ref{eqn:duffing-AB})--(\ref{eqn:data}) satisfies (\ref{th:main-alternative}) with $\sigma=0$ is a nonempty closed subset of $\dath$.
	
\end{enumerate}

\end{thm}

When there is no external force, or the external force vanishes in the limit, then all solutions tend to one of the three stationary points of the functional $\mathcal{E}_{A,B}(u)$.

\begin{cor}[Asymptotically unforced equation]\label{cor}

Let us consider problem (\ref{eqn:duffing-AB})--(\ref{eqn:data}) under the standard assumptions presented before Theorem~\ref{thm:main}. 

Let us assume in addition that
\begin{equation}
\lim_{t\to+\infty}f(t)=0
\quad\quad
\mbox{in }H.
\nonumber
\end{equation}

Then there exists $\sigma\in\{-\sigma_{0},0,\sigma_{0}\}$ such that
\begin{equation}
\lim_{t\to +\infty}\left(|u'(t)|+|B(u(t)-\sigma e_{1}|\strut\right)=0.
\nonumber
\end{equation}

\end{cor}

\paragraph{\textmd{\textit{Application to the beam equation with clamped ends}}} 

The abstract results stated in Theorem~\ref{thm:main} and Corollary~\ref{cor} can be applied to the beam equation (\ref{eqn:duffing-concr}) with clamped ends. To this end, it is enough to show that the relevant operators $A$ and $B$ fit in the framework of Definition~\ref{defn:gap}. 

\begin{prop}[Abstract setting for the beam with clamped ends]\label{prop:concrete}

Let us consider the Hilbert space $H=L^{2}((0,1))$, the self-adjoint positive unbounded operator $A$ on $H$ defined by $Au=-u_{xx}$ with domain
\begin{equation}
D(A):=H^{2}((0,1))\cap H^{1}_{0}((0,1)),
\nonumber
\end{equation}
and the self-adjoint positive unbounded operator $B$ on $H$ such that $B^{2}u=u_{xxxx}$ with domain
\begin{equation}
D(B^{2}):=H^{4}((0,1))\cap H^{2}_{0}((0,1)).
\nonumber
\end{equation}

Then it turns out that $(A,B)\in\gaph$ with $\lambda_{1}=4\pi^{2}$ and $\lambda_{2}=4\alpha_{1}^{2}$, where $\alpha_{1}\in(\pi,3\pi/2)$ is the smallest positive real solution to the equation $\tan\alpha=\alpha$.

\end{prop}


\setcounter{equation}{0}
\section{Auxiliary results}\label{sec:proof-sketch}

In this section we state four auxiliary results that correspond to the key steps in the proof of Theorem~\ref{thm:main}.

To begin with, we introduce the operator
\begin{equation}
Ru:=u-2\langle u,e_{0}\rangle e_{0}
\qquad
\forall u\in H,
\label{defn:Ru}
\end{equation}
where $e_{0}$ is a unit vector satisfying (\ref{defn:e0}). We observe that $R$ is the bounded operator on $H$ such that $Re_{0}=-e_{0}$ and $Rv=v$ for every $v$ orthogonal to $e_{0}$. Then we set
\begin{equation}
\delta:=\frac{1}{2}\left(1+\frac{2}{\mu_{2}^{1/2}}|A^{1/2}e_{0}|\right)^{-1},
\label{defn:delta}
\end{equation}
where $\mu_{2}$ is the constant that appears in (\ref{defn:mu-2}), and we define the operator
\begin{equation}
Pu:=u+\delta Ru
\qquad
\forall u\in H.
\label{defn:Pu}
\end{equation}

Finally, we choose a positive real number $\gamma_{0}$ such that
\begin{equation}
\gamma_{0}\leq\frac{1}{2(5+2\delta)},
\qquad\qquad
\gamma_{0}\leq\frac{\delta}{(5+\delta)(1+\delta)}\min\left\{\mu_{2}^{2}\mu_{1}\mu_{3},\lambda_{0}\right\},
\label{defn:gamma0}
\end{equation}
where $\mu_{1}$, $\lambda_{0}$ and $\mu_{3}$ are the constants that appear in (\ref{defn:mu-1}), (\ref{defn:e0}) and (\ref{defn:e0+}), respectively.

For every solution $u(t)$ to (\ref{eqn:duffing-AB}), we consider the classical energy $E(t)$ defined by (\ref{defn:E}), and the modified energy
\begin{equation}
F(t):=E(t)+2\gamma_{0}\langle Pu(t),u'(t)\rangle+\gamma_{0}\langle Pu(t),u(t)\rangle.
\label{defn:F}
\end{equation}

In the first result we prove that the energy $F(t)$ of solutions to (\ref{eqn:duffing-AB}) is bounded for $t$ large in terms of the norm of the forcing term. As a consequence, all solutions are bounded in $\dath$.

\begin{prop}[Ultimate bound on solutions]\label{prop:ultimate}

Let us consider problem (\ref{eqn:duffing-AB})--(\ref{eqn:data}) under the standard assumptions presented before Theorem~\ref{thm:main}. Let $F(t)$ be the energy defined in (\ref{defn:F}).

Then there exists a positive constant $M_{1}$, independent of the forcing term $f(t)$ and of the solution $u(t)$, such that
\begin{equation}
\limsup_{t\to+\infty}F(t)\leq M_{1}\limsup_{t\to+\infty}|f(t)|^{2},	
\label{th:ultimate}
\end{equation}
and there exist two positive constants $M_{2}$ and $M_{3}$, again independent of $f(t)$ and $u(t)$, such that
\begin{equation}
\limsup_{t\to+\infty}\left(|u'(t)|^{2}+|Bu(t)|^{2}\right)\leq M_{2}+M_{3}\limsup_{t\to+\infty}|f(t)|^{2}.	
\label{th:ultimate-u}
\end{equation}
\end{prop}

In the second result we deal with solutions $u(t)$ such that $|u(t)|$ is eventually smaller than a universal constant. We show that the norm of these solutions in the energy space is asymptotically bounded by the norm in $H$ of the forcing term.

\begin{prop}[Solutions in the unstable regime]\label{prop:unstable}

Let us consider problem (\ref{eqn:duffing-AB})--(\ref{eqn:data}) under the standard assumptions presented before Theorem~\ref{thm:main}. 

Then there exist two positive constants $\beta_{0}$ and $M_{4}$, independent of the forcing term $f(t)$ and of the solution $u(t)$, for which the following implication is true:
\begin{equation}
\boxa{\begin{array}{c}
\displaystyle\limsup_{t\to+\infty}|f(t)|\leq 1 \\
\noalign{\vspace{1ex}}
\displaystyle\limsup_{t\to+\infty}|u(t)|\leq\beta_{0}
\end{array}}
\ \Longrightarrow\ 
\boxa{\limsup_{t\to+\infty}\left(|u'(t)|+|Bu(t)|\strut\right)\leq M_{4}\limsup_{t\to+\infty}|f(t)|}
\label{th:unstable}
\end{equation}

\end{prop}

In the third result we deal with solutions that, \emph{at a given time}, are close to one of the stable stationary points of the functional (\ref{defn:funct-AB}). Here ``close to a stationary point'' means that the energy $E(t)$ is negative at the given time. We show that these solutions lie eventually in a neighborhood of the same stationary point, within a distance depending on the norm in $H$ of the forcing term.

\begin{prop}[Solutions in the stable regime]\label{prop:stable}

Let us consider problem (\ref{eqn:duffing-AB})--(\ref{eqn:data}) under the standard assumptions presented before Theorem~\ref{thm:main}.  Let $e_{1}$ be the vector that appears in (\ref{hp:e1}) and (\ref{hp:e1-perp}), and let $E(t)$ be the energy defined in (\ref{defn:E}).

Then for every $\eta>0$ there exist two constants $\ep_{1}>0$ and $M_{5}>0$, independent of the forcing term $f(t)$ and of the solution $u(t)$, for which the following implication is true: 
\begin{equation}
\begin{array}{c}
\fbox{$\begin{array}{c}
\exists\, T_{0}\geq 0 \mbox{ such that } \\[1ex]
\displaystyle\sup_{t\geq T_{0}}|f(t)|\leq\ep_{1},\ E(T_{0})<-\eta
\end{array}$}   \\[4ex]
\Downarrow   \\[1ex]
\fbox{$\begin{array}{c}
\exists\, \sigma\in\{\sigma_{0},-\sigma_{0}\} \mbox{ such that } \\[1ex]
\displaystyle\limsup_{t\to+\infty}\left(|u'(t)|+|B(u(t)-\sigma e_{1})|\strut\right)\leq M_{5}\limsup_{t\to+\infty}|f(t)|
\end{array}$}
\end{array}
\label{th:stable}
\end{equation}

Moreover, when the assumptions in the upper box of (\ref{th:stable}) are satisfied, it turns out that
\begin{equation}
\langle u(t),Ae_{1}\rangle\neq 0
\qquad
\forall t\geq T_{0},
\label{th:stable-sign}
\end{equation}
and as a consequence the sign of $\sigma$ coincides with the sign of $\langle u(T_{0}),Ae_{1}\rangle$.

\end{prop}

In the last result we show that any two solutions to (\ref{eqn:duffing-AB}) that are close enough to the same stationary point of the functional (\ref{defn:funct-AB}) are actually asymptotic to each other.

\begin{prop}[Close solutions are asymptotic to each other]\label{prop:asymptotic}

Let us consider problem (\ref{eqn:duffing-AB})--(\ref{eqn:data}) under the standard assumptions presented before Theorem~\ref{thm:main}. 

Then there exists $r_{0}>0$ with the following property: if $u(t)$ and $v(t)$ are two solutions to (\ref{eqn:duffing-AB}), and there exists $\sigma\in\{-\sigma_{0},0,\sigma_{0}\}$ such that
\begin{equation}
\limsup_{t\to +\infty}\left(|u'(t)|+|B(u(t)-\sigma e_{1})|+|v'(t)|+|B(v(t)-\sigma e_{1})|\strut\right)\leq r_{0},
\label{hp:asymptotic}
\end{equation}
then $u(t)$ and $v(t)$ are asymptotic to each other in the sense of (\ref{th:main-asymptotic}).

\end{prop}


\setcounter{equation}{0}
\section{Some linear algebra in infinite dimensions}\label{sec:linear-algebra}

A fundamental tool in~\cite{GGH:duffing} was considering the components of a solution $u(t)$ with respect to the eigenspaces of $A$. Due to the presence of two operators, in this paper we are forced to consider different decompositions of the Hilbert space $H$ in different parts of the proof. In this section we introduce the decompositions that we need in the sequel, we state their basic properties, and we prove the linear algebra results of Proposition~\ref{prop:gaph}.

\subsection{Decomposition of the space in the stable regime}

Let $e_{1}$ be the vector that appears in (\ref{hp:e1}) and (\ref{hp:e1-perp}). We consider the decomposition
\begin{equation}
H:=W_{-}\oplus W_{+},
\label{defn:W+-}
\end{equation}
where $W_{-}$ is the one-dimensional subspace spanned by $e_{1}$, and $W_{+}$ is the subspace orthogonal to $Ae_{1}$. We point out that this is a direct sum,  in general not orthogonal in the sense of $H$ (but orthogonality is true in the sense of $D(A^{1/2})$ for vectors belonging to $D(A^{1/2})$, both projections remaining in $D(A^{1/2})$). Every vector $u\in H$ can be written in a unique way in the form 
\begin{equation}
u=\alpha e_{1}+w,
\nonumber
\end{equation}
where $\alpha\in\re$ and $w\in W_{+}$ are given by
\begin{equation}
\alpha:=\langle u,Ae_{1}\rangle
\qquad\mbox{and}\qquad
w:=u-\alpha e_{1}.
\nonumber
\end{equation}

Due to (\ref{hp:e1}), it turns out that
\begin{equation}
\langle w,B^{2}e_{1}\rangle=\lambda_{1}\langle w,Ae_{1}\rangle=0
\qquad
\forall w\in W_{+},
\label{eqn:w-e1}
\end{equation}
from which it follows that 
\begin{equation}
|A^{1/2}u|^{2}=\alpha^{2}+|A^{1/2}w|^{2}
\label{A=a+w}
\end{equation}
and
\begin{equation}
|Bu|^{2}=\alpha^{2}|Be_{1}|^{2}+|Bw|^{2}=\lambda_{1}\alpha^{2}+|Bw|^{2}
\label{B=a+w}
\end{equation}
for every $u\in D(B)$. Moreover, from (\ref{hp:e1-perp}) we deduce that
\begin{equation}
|Bw|^{2}=\frac{\lambda_{2}-\lambda}{\lambda_{2}}|Bw|^{2}+\frac{\lambda}{\lambda_{2}}|Bw|^{2}\geq\frac{\lambda_{2}-\lambda}{\lambda_{2}}|Bw|^{2}+\lambda|A^{1/2}w|^{2}
\nonumber
\end{equation}
for every $w\in D(B)\cap W_{+}$, and in particular 
\begin{equation}
|Bw|^{2}-\lambda|A^{1/2}w|^{2}\geq\frac{\lambda_{2}-\lambda}{\lambda_{2}}|Bw|^{2}
\qquad
\forall w\in D(B)\cap W_{+}.
\label{est:Bw}
\end{equation}

From (\ref{A=a+w}), (\ref{B=a+w}), and (\ref{defn:sigma-0}) we deduce that, if $\lambda>\lambda_{1}$,
\begin{eqnarray}
|Bu|^{2}-\lambda|A^{1/2}u|^{2} & = & (\lambda_{1}-\lambda)\alpha^{2}+|Bw|^{2}-\lambda|A^{1/2}w|^{2} 
\nonumber  \\[1ex]
& = & -\sigma_{0}^{2}\alpha^{2}+|Bw|^{2}-\lambda|A^{1/2}w|^{2}.
\label{eqn:Bu-Au}
\end{eqnarray}

\subsection{Proof of Proposition~\ref{prop:gaph}}

\paragraph{\textmd{\textit{Statement~(1)}}}

We have to prove that the operator $B^{2}-\lambda A$ is self-adjoint with domain $D(B^{2})$. To this end, it is enough to prove that the operator $B^{2}u-\lambda Au+mu$ is symmetric and maximal monotone with domain $D(B^{2})$ for some real number $m$ (here we exploit that a symmetric maximal monotone operator is self-adjoint, and that the sum of a self-adjoint operator and a bounded symmetric operator is again self-adjoint).

The symmetry is trivial, and therefore we can limit ourselves to check monotonicity and maximality.
\begin{itemize}

\item  We claim that $B^{2}u-\lambda Au+mu$ is monotone when $m$ is large enough. If $\lambda\leq 0$ the conclusion is true even with $m=0$. If $\lambda>0$ we exploit the inequality
$$|A^{1/2}u|^{2}=\langle Au,u\rangle\leq\frac{\ep}{2}|Au|^{2}+\frac{1}{2\ep}|u|^{2}
\qquad
\forall u\in D(A)\quad\forall\ep>0,$$
and from (\ref{defn:mu-1}) we deduce that
\begin{eqnarray*}
\langle B^{2}u-\lambda Au+mu,u\rangle & = & |Bu|^{2}-\lambda|A^{1/2}u|^{2}+m|u|^{2} \\[0.5ex]
& \geq & \mu_{1}|Au|^{2}-\frac{\lambda\ep}{2}|Au|^{2}-\frac{\lambda}{2\ep}|u|^{2}+m|u|^{2}.
\end{eqnarray*}

At this point it is enough to choose $\ep=2\mu_{1}\lambda^{-1}$ and $\displaystyle m\geq \frac{\lambda}{2\ep}$.

\item  We claim that $B^{2}u-\lambda Au+mu$ is surjective from $D(B^{2})$ to $H$ when $m$ is large enough, namely that for every $f\in H$ there exists $u\in D(B^{2})$ such that
\begin{equation}
B^{2}u-\lambda Au+mu=f.
\label{eqn:pre-fp}
\end{equation}

To this end, we exploit a fixed point technique. For every $v\in D(A)$, the equation
\begin{equation}
B^{2}u+mu=\lambda Av+f
\label{eqn:fixed-point}
\end{equation}
has a unique solution $u=T(v)\in D(B^{2})\subseteq D(A)$. This defines a function $T:D(A)\to D(A)$. Any fixed point of $T$ lies in $D(B^{2})$, and is a solution to (\ref{eqn:pre-fp}). Let $v_{1}$ and $v_{2}$ be in $D(A)$, and let us set $u_{1}=T(v_{1})$ and $u_{2}=T(v_{2})$. From (\ref{eqn:fixed-point}) we obtain that
\begin{eqnarray*}
|B(u_{2}-u_{1})|^{2}+m|u_{2}-u_{1}|^{2}  & =  & \lambda\langle A(v_{2}-v_{1}),u_{2}-u_{1}\rangle \\[1ex]
& \leq & \lambda|A(v_{2}-v_{1})|\cdot|u_{2}-u_{1}| \\[0.5ex]
 &  \leq  & \frac{\lambda^{2}}{2m}|A(v_{2}-v_{1})|^{2}+\frac{m}{2}|u_{2}-u_{1}|^{2}. 
\end{eqnarray*}

Recalling (\ref{defn:mu-1}), when $m$ is large enough we deduce that
$$|A(u_{2}-u_{1})|^{2}+|u_{2}-u_{1}|^{2}\leq\frac{1}{2}|A(v_{2}-v_{1})|^{2},$$
and hence $T$ is a contraction in the Hilbert space $D(A)$. This proves that the operator $B^{2}u-\lambda Au+mu$ is also maximal with domain $D(B^{2})$.

\end{itemize}

\begin{rmk}
\begin{em}
The sum of a maximal monotone operator $L_{1}$ and a monotone operator $L_{2}$ dominated by $L_{1}$, with the same (or larger) domain, is again maximal monotone (see~\cite[Proposition~2.10]{brezis:OMM}). Therefore, we can give an alternative proof by writing
$$B^{2}u-\lambda Au+mu=\frac{2}{3}B^{2}u+\frac{1}{3}B^{2}u-\lambda Au+mu,$$
and applying the abstract result with $L_{1}u:=\frac{2}{3}B^{2}u$ and $L_{2}u:=\frac{1}{3}B^{2}u-\lambda Au+mu$. It is enough to check that $L_{2}$ is monotone and dominated by $L_{1}$, and this can be done as we did in the first item of the previous proof.
\end{em}
\end{rmk}

\paragraph{\textmd{\textit{Statement~(2)}}}

Let us write any element $u\in H$ in the form $u=\alpha e_{1}+w$ according to the direct sum (\ref{defn:W+-}). From (\ref{A=a+w}), (\ref{B=a+w}) and (\ref{hp:e1-perp}) it follows that
\begin{eqnarray*}
|Bu|^{2}-\lambda|A^{1/2}u|^{2} & = & (\lambda_{1}-\lambda)\alpha^{2}+|Bw|^{2}-\lambda|A^{1/2}w|^{2} \\[0.5ex]
 & \geq & (\lambda_{1}-\lambda)\alpha^{2}+(\lambda_{2}-\lambda)|A^{1/2}w|^{2}  \\[0.5ex]
 & \geq & 0,
\end{eqnarray*}
with strict inequality if either $\alpha\neq 0$ or $w\neq 0$. This proves (\ref{th:prop-sub}).

\paragraph{\textmd{\textit{Statement~(3) -- Computation of the negative inertia index}}}

For every $\lambda>\lambda_{1}$, the operator $B^{2}-\lambda A$ is negative in the one-dimensional subspace of $H$ generated by $e_{1}$. 

We claim that, if $\lambda<\lambda_{2}$, the same operator cannot be negative, or even just less than or equal to~0, in any subspace of $H$ of dimension at least two. Indeed, any such subspace contains a vector $v\neq 0$ with $\langle v,Ae_{1}\rangle=0$, and for this vector it turns out that
\begin{equation}
|Bv|^{2}-\lambda|A^{1/2}v|^{2}>|Bv|^{2}-\lambda_{2}|A^{1/2}v|^{2}\geq 0,
\nonumber
\end{equation}
where the second inequality follows from (\ref{hp:e1-perp}).

\paragraph{\textmd{\textit{Statement~(3) -- Existence of an eigenvector with negative eigenvalue}}}

According to the spectral theory (see for example~\cite[Theorem VIII.4]{reed-simon}), we can identify $H$ with $L^{2}(M,\mu)$ for some measure space $(M,\mu)$ in such a way that under this identification the operator $B^{2}-\lambda A$ becomes a multiplication operator. This means that there exists a measurable function $\lambda(\xi)$ in $M$ with the property that, if $\widehat{u}(\xi)\in L^{2}(M,\mu)$ corresponds to some $u\in H$, then $\lambda(\xi)\widehat{u}(\xi)$ corresponds to $B^{2}u-\lambda Au$.

Let us consider the set
\begin{equation}
N:=\{\xi\in M:\lambda(\xi)<0\}.
\nonumber
\end{equation}

We claim that $\mu(N)>0$, and $\lambda(\xi)$ is equal to some negative constant $-\lambda_{0}$ for almost every $\xi\in N$. If we prove these claims, then the vector $e_{0}\in H$ that under the identification corresponds to the characteristic function of $N$ is an eigenvector of $B^{2}-\lambda A$ with  eigenvalue $-\lambda_{0}$.

In order to prove that $\mu(N)>0$ it is enough to observe that otherwise the operator $B^{2}-\lambda A$ would be nonnegative in $H$.

In order to prove that $\lambda(\xi)$ is essentially constant in $N$, let us assume that this is not the case. Then there exists a real number $\lambda_{*}<0$ such that the two sets
\begin{equation}
N_{1}:=\{\xi\in N:\lambda(\xi)\leq\lambda_{*}\}
\qquad\mbox{and}\qquad
N_{2}:=\{\xi\in N:\lambda(\xi)>\lambda_{*}\}
\nonumber
\end{equation}
have positive measure. In this case the two vectors $u_{1}$ and $u_{2}$ corresponding to the characteristic functions of $N_{1}$ and $N_{2}$ would be two orthogonal vectors that span a two-dimensional subspace of $H$ where $B^{2}-\lambda A$ is negative, and we already know that this is not possible when $\lambda\in(\lambda_{1},\lambda_{2})$.

\paragraph{\textmd{\textit{Statement~(3) -- Estimate in the orthogonal space}}}

Let us prove that (\ref{defn:e0+}) holds true if we choose
\begin{equation}
\mu_{3}<\min\left\{\frac{\lambda_{2}-\lambda}{\lambda_{2}+\lambda},\frac{\lambda_{0}}{2\lambda|A^{1/2}e_{0}|^{2}+\lambda_{0}}\right\}.
\label{defn:mu3}
\end{equation}

To this end, let us assume that it is not the case. Then there exists $v\in D(B)$ with
\begin{equation}
\langle v,e_{0}\rangle=0
\qquad\mbox{and}\qquad
|Bv|^{2}-\lambda|A^{1/2}v|<\mu_{3}|Bv|^{2}.
\label{defn:mu3-bis}
\end{equation}

Let us set
\begin{equation}
\widehat{\lambda}:=\frac{1+\mu_{3}}{1-\mu_{3}}\lambda,
\label{defn:hat-lambda}
\end{equation}
and let us observe that $\widehat{\lambda}<\lambda_{2}$ because of the first request in the definition of $\mu_{3}$. 

Now we show that the operator $B^{2}-\widehat{\lambda}A$ is less than or equal to~0 on the two-dimensional subspace of $H$ spanned by $e_{0}$ and $v$, which we already shown to be absurd. To this end, we take a generic vector $u=\alpha e_{0}+\beta v$, and with some computations we obtain that
\begin{eqnarray*}
|Bu|^{2}-\widehat{\lambda}|A^{1/2}u|^{2} & = & |Bu|^{2}-\lambda|A^{1/2}u|^{2}-(\widehat{\lambda}-\lambda)|A^{1/2}u|^{2} \\[1ex]
& = & -\lambda_{0}\alpha^{2}+\beta^{2}\left(|Bv|^{2}-\lambda|A^{1/2}v|^{2}\right) \\[1ex]
&& -(\widehat{\lambda}-\lambda)\left(\alpha^{2}|A^{1/2}e_{0}|^{2}+\beta^{2}|A^{1/2}v|^{2}+2\alpha\beta\langle A^{1/2}e_{0},A^{1/2}v\rangle\right).
\end{eqnarray*}

Now from (\ref{defn:mu3}) and (\ref{defn:mu3-bis}) we deduce that
\begin{equation}
|Bv|^{2}<\frac{\lambda}{1-\mu_{3}}|A^{1/2}v|^{2}
\qquad\mbox{and}\qquad
\frac{\lambda_{0}}{|A^{1/2}e_{0}|^{2}}+(\widehat{\lambda}-\lambda)\geq\frac{4\lambda\mu_{3}}{1-\mu_{3}},
\nonumber
\end{equation}
and keeping (\ref{defn:hat-lambda}) into account we conclude that
\begin{eqnarray*}
|Bu|^{2}-\widehat{\lambda}|A^{1/2}u|^{2} & \leq & -\alpha^{2}\left\{\frac{\lambda_{0}}{|A^{1/2}e_{0}|^{2}}+(\widehat{\lambda}-\lambda)\right\}|A^{1/2}e_{0}|^{2}-\beta^{2}\left\{\widehat{\lambda}-\frac{\lambda}{1-\mu_{3}}\right\}|A^{1/2}v|^{2} \\[1ex]
&& \mbox{}+2(\widehat{\lambda}-\lambda)|\alpha|\cdot|\beta|\cdot|A^{1/2}e_{0}|\cdot|A^{1/2}v| \\[1ex]
& \leq & -\frac{\lambda\mu_{3}}{1-\mu_{3}}\left(2|\alpha|\cdot|A^{1/2}e_{0}|-|\beta|\cdot|A^{1/2}v|\right)^{2} \\[0.5ex]
& \leq & 0.
\end{eqnarray*}

\paragraph{\textmd{\textit{Statement~(4)}}}

Let us set $\mu:=|A^{1/2}u|^{2}$, and let us look for nonzero solutions to equation
\begin{equation}
B^{2}u-(\lambda-\mu)Au=0.
\label{eqn:B-mu-lambda}
\end{equation}

Let us write as usual $u=\alpha e_{1}+w$ according to the direct sum (\ref{defn:W+-}). Then equation (\ref{eqn:B-mu-lambda}) reduces to 
\begin{equation}
\alpha B^{2}e_{1}-(\lambda-\mu)\alpha Ae_{1}+B^{2}w-(\lambda-\mu)Aw=0.
\nonumber
\end{equation}

Due to (\ref{hp:e1}), taking the scalar product of this equation with $w$ we obtain that
$$|Bw|^{2}=(\lambda-\mu)|A^{1/2}w|^{2}.$$ 

Since $\lambda-\mu<\lambda_{2}$, this is impossible if $w\neq 0$ because of (\ref{hp:e1-perp}). It follows that $u=\alpha e_{1}$ for some $\alpha\neq 0$, and 
\begin{equation}
\alpha B^{2}e_{1}-(\lambda-\mu)\alpha Ae_{1}=0.
\nonumber
\end{equation}

Keeping (\ref{hp:e1}) into account, we conclude that
\begin{equation}
\lambda_{1}=\lambda-\mu=\lambda-|A^{1/2}u|^{2}=\lambda-\alpha^{2}{},
\nonumber
\end{equation}
which implies that $\alpha=\pm\sigma_{0}$, with $\sigma_{0}$ given by (\ref{defn:sigma-0}). 
 
\paragraph{\textmd{\textit{Statement~(5)}}}

If $\lambda_{2}$ is the largest constant for which (\ref{hp:e1-perp}) holds true, then for every $\lambda>\lambda_{2}$ there exists a vector $v\neq 0$ (possibly depending on $\lambda$) such that $\langle v,Ae_{1}\rangle=0$ and $|Bv|^{2}<\lambda|A^{1/2}v|^{2}$. At this point it turns out that (\ref{th:prop-4}) holds true in the two-dimensional subspace of $H$ spanned by $e_{1}$ and $v$.

\subsection{Decomposition of the space in the unstable regime}

Let us consider the operator $L:=B^{2}-\lambda A$, which we know to be self-adjoint with domain $D(L)=D(B^{2})$. Since $\lambda\in(\lambda_{1},\lambda_{2})$, from statement~(3) of Proposition~\ref{prop:gaph} we know that $L$ has a negative eigenvalue $-\lambda_{0}$. Given a corresponding eigenvector $e_{0}$ with unit norm, we write $H$ as a direct orthogonal sum
\begin{equation}
H:=H_{-}\oplus H_{+},
\label{defn:H+-}
\end{equation} 
where $H_{-}$ is the one-dimensional subspace spanned by $e_{0}$, and $H_{+}$ is the subspace orthogonal to $e_{0}$. In this way any vector $u\in H$ can be written in a unique way as the sum of a low-frequency component $u_{-}\in H_{-}$ and a high-frequency component $u_{+}\in H_{+}$, where of course
\begin{equation}
u_{-}:=\langle u,e_{0}\rangle\, e_{0}
\qquad\mbox{and}\qquad
u_{+}:=u-u_{-}.
\nonumber
\end{equation}

We point out that $H_{-}$ and $H_{+}$ are invariant subspaces for $L$, but they are not necessarily invariant spaces for $A$ or $B$. 

From (\ref{defn:e0+}) we know that $L$ is a positive operator when restricted to $H_{+}$, and
\begin{equation}
|L^{1/2}u|^{2}\geq\mu_{3}|Bu|^{2}
\qquad
\forall u\in D(B)\cap H_{+}.
\label{est:L-B-A}
\end{equation}

Since from (\ref{defn:mu-1}) and (\ref{defn:mu-2}) we know that
\begin{equation}
|Bu|^{2}\geq\mu_{1}|Au|^{2}\geq\mu_{1}\mu_{2}|A^{1/2}u|^{2}\geq\mu_{1}\mu_{2}^{2}|u|^{2}
\qquad
\forall u\in D(B),
\label{est:B2A}
\end{equation}
from (\ref{est:L-B-A}) we can derive estimates for $|u|$, $|A^{1/2}u|$ and $|Au|$ in terms of $|L^{1/2}u|$ for every $u\in D(B)\cap H_{+}$.


\setcounter{equation}{0}
\section{Proof of auxiliary results}\label{sec:proof-prop}

\subsection{Useful ultimate bounds}

In this subsection we recall three results concerning ultimate bounds that are crucial in the sequel. For a proof we refer to~\cite[Section~4.1]{GGH:duffing} and to the references quoted therein.

\begin{lemma}\label{lemma:limsup}

For every bounded function $\varphi:[0,+\infty)\to\re$ of class $C^{1}$ there exists a sequence $t_{n}\to +\infty$ of nonnegative real numbers such that
$$\lim_{n\to+\infty}\varphi(t_{n})=\limsup_{t\to+\infty}\varphi(t)
\quad\quad\mbox{and}\quad\quad
\lim_{n\to+\infty}\varphi'(t_{n})=0.$$

\end{lemma}  

\begin{lemma}\label{lemma:ODE-limsup}

Let $m$ be a positive real number, let $\psi:[0,+\infty)\to\re$ be a continuous function, and let $y\in C^2([0,+\infty), \re) $ be a solution to
\begin{equation}
y''(t)+y'(t)-my(t)=\psi(t).
\nonumber
\end{equation}

Let us assume that both $\psi(t)$ and $y(t)$ are bounded.

Then it turns out that
\begin{equation}
\limsup_{t\to+\infty}|y(t)|\leq\frac{1}{m}\limsup_{t\to+\infty}|\psi(t)|,
\nonumber
\end{equation}
\begin{equation}
\limsup_{t\to+\infty}|y'(t)|\leq 2\limsup_{t\to+\infty}|\psi(t)|.
\nonumber
\end{equation}

\end{lemma}

\begin{lemma}\label{lemma:PDE-limsup}

Let $X$ be a Hilbert space, and let $\mathcal{L}$ be a self-adjoint linear operator on $X$ with dense domain $D(\mathcal{L})$. Let us assume that there exists a constant $m>0$ such that
\begin{equation}
\langle \mathcal{L} x,x\rangle\geq m|x|^{2}
\quad\quad
\forall x\in D(\mathcal{L}).
\nonumber
\end{equation}

Let $\psi:[0,+\infty)\to X$ be a bounded continuous function, and let 
$$y\in C^{0}\left([0,+\infty),D(\mathcal{L}^{1/2})\right)\cap C^{1}\left([0,+\infty),X\right)$$ 
be a solution to
\begin{equation}
y''(t)+y'(t)+\mathcal{L} y(t)=\psi(t).
\nonumber
\end{equation}

Then it turns out that
\begin{equation}
\limsup_{t\to+\infty}\left(|y'(t)|^{2}+|\mathcal{L}^{1/2}y(t)|^{2}\right)\leq 9\max\left\{1,\frac{1}{m}\right\}\cdot\limsup_{t\to+\infty}|\psi(t)|^2.
\nonumber
\end{equation}

\end{lemma}


\subsection{Useful estimates for functionals and energies}

In this subsection we collect some identities and inequalities that are needed several times in the sequel.

Let $R$ and $P$ be the operators defined in (\ref{defn:Ru}) and (\ref{defn:Pu}). Let us consider the orthogonal direct sum (\ref{defn:H+-}). Then it turns out that
\begin{equation}
Ru=\left\{
\begin{array}{cc}
-u & \mbox{if }u\in H_{-}, \\[0.5ex]
u & \mbox{if }u\in H_{+},
\end{array}
\right.
\qquad\mbox{and}\qquad
Pu=\left\{
\begin{array}{cc}
(1-\delta)u & \mbox{if }u\in H_{-}, \\[0.5ex]
(1+\delta)u & \mbox{if }u\in H_{+}.
\end{array}
\right.
\nonumber
\end{equation}

As a consequence we obtain that
\begin{equation}
(1-\delta)|u|^{2}\leq\langle Pu,u\rangle\leq(1+\delta)|u|^{2}
\qquad
\forall u\in H,
\label{est:pu,u-u,u}
\end{equation}
and
\begin{equation}
(1-\delta)\langle Pu,u\rangle\leq|Pu|^{2}\leq(1+\delta)\langle Pu,u\rangle
\qquad
\forall u\in H.
\label{est:Pu,Pu-Pu,u}
\end{equation}

Let us consider the functional $\mathcal{E}_{A,B}(u)$ defined in (\ref{defn:funct-AB}). Since 
$$\frac{1}{4}x^{4}-\frac{\lambda}{2}x^{2}\geq -\frac{\lambda^{2}}{4}
\qquad
\forall x\in\re,$$
we obtain that
\begin{equation}
\mathcal{E}_{A,B}(u)\geq\frac{1}{2}|Bu|^{2}-\frac{\lambda^{2}}{4}
\qquad
\forall u\in D(B),
\label{est:EAB}
\end{equation}
and analogously
\begin{equation}
E(t)\geq\frac{1}{2}|u'(t)|^{2}+\frac{1}{2}|Bu(t)|^{2}-\frac{\lambda^{2}}{4}.
\label{est:E-below}
\end{equation}

Let us write now $u$ in the form $\alpha e_{1}+w$ according to the direct sum (\ref{defn:W+-}). Then from (\ref{A=a+w}), (\ref{eqn:Bu-Au}), and (\ref{est:Bw}) we deduce that
\begin{eqnarray}
\mathcal{E}_{A,B}(u) & = & \frac{1}{4}\alpha^{4}-\frac{1}{2}\sigma_{0}^{2}\alpha^{2}+\frac{1}{2}|Bw|^{2}-\frac{\lambda}{2}|A^{1/2}w|^{2} 
\nonumber  \\[1ex]
& & \mbox{}+\frac{1}{4}|A^{1/2}w|^{4}+\frac{1}{2}\alpha^{2}{}\cdot|A^{1/2}w|^{2} 
\nonumber   \\[1ex]
& \geq & \frac{1}{4}\alpha^{4}-\frac{1}{2}\sigma_{0}^{2}\alpha^{2}+\frac{1}{2}\frac{\lambda_{2}-\lambda}{\lambda_{2}}|Bw|^{2}.
\label{est:E-below-aw}
\end{eqnarray}

This shows in particular that
\begin{equation}
\mathcal{E}_{A,B}(u)\geq -\frac{1}{4}\sigma_{0}^{4}{}
\qquad
\forall u\in D(B),
\nonumber
\end{equation}
with equality if and only if $u=\pm\sigma_{0}e_{1}$, and in addition
\begin{equation}
\mathcal{E}_{A,B}(u)+\frac{1}{4}\sigma_{0}^{4}{}\geq\frac{1}{4}\left(\alpha^{2}-\sigma_{0}^{2}\right)^{2}{}+\frac{1}{2}\frac{\lambda_{2}-\lambda}{\lambda_{2}}|Bw|^{2}.
\label{eqn:Bw}
\end{equation}


\subsection{Proof of Proposition~\ref{prop:ultimate}}

\paragraph{\textmd{\textit{Estimates on $F(t)$ from above and below}}}

Let us consider the energy $F(t)$ defined in (\ref{defn:F}). We prove that
\begin{equation}
F(t)\leq\frac{3}{4}|u'(t)|^{2}+\frac{1}{2}|Bu(t)|^{2}-\frac{\lambda}{2}|A^{1/2}u(t)|^{2}+\frac{1}{4}|A^{1/2}u(t)|^{4}+2\gamma_{0}(1+\delta)|u(t)|^{2}
\label{est:F-above}
\end{equation}
and
\begin{equation}
F(t)\geq\frac{1}{4}|u'(t)|^{2}+\frac{1}{2}|Bu(t)|^{2}-\frac{\lambda}{2}|A^{1/2}u(t)|^{2}+\frac{1}{4}|A^{1/2}u(t)|^{4}.
\label{est:F-below}
\end{equation}

Indeed, from (\ref{est:Pu,Pu-Pu,u}) it follows that
\begin{equation}
|2\langle Pu(t),u'(t)\rangle|\leq(1+\delta)|u'(t)|^{2}+\frac{1}{1+\delta}|Pu(t)|^{2}\leq(1+\delta)|u'(t)|^{2}+\langle Pu(t),u(t)\rangle.
\nonumber
\end{equation}

Plugging this inequality into (\ref{defn:F}) we deduce that
\begin{eqnarray}
F(t) & \leq & \left(\frac{1}{2}+\gamma_{0}(1+\delta)\right)|u'(t)|^{2}+\frac{1}{2}|Bu(t)|^{2}-\frac{\lambda}{2}|A^{1/2}u(t)|^{2} 
\nonumber  \\[0.5ex]
& & \mbox{}+\frac{1}{4}|A^{1/2}u(t)|^{4}+2\gamma_{0}\langle Pu(t),u(t)\rangle
\label{est:F-above-pre}
\end{eqnarray}
and
\begin{equation}
F(t)\geq\left(\frac{1}{2}-\gamma_{0}(1+\delta)\right)|u'(t)|^{2}+\frac{1}{2}|Bu(t)|^{2}-\frac{\lambda}{2}|A^{1/2}u(t)|^{2}+\frac{1}{4}|A^{1/2}u(t)|^{4}.
\label{est:F-below-pre}
\end{equation}

Finally, we observe that (\ref{defn:delta}) and (\ref{defn:gamma0}) imply in particular that $\gamma_{0}(1+\delta)\leq 1/4$. At this point, (\ref{est:F-above}) follows from (\ref{est:F-above-pre}) and (\ref{est:pu,u-u,u}), while  (\ref{est:F-below}) follows from (\ref{est:F-below-pre}).

\paragraph{\textmd{\textit{Estimates for the operator $R$}}}

We show that
\begin{equation}
\delta|\langle Ru,Au\rangle|\leq\frac{1}{2}|A^{1/2}u|^{2}
\qquad
\forall u\in D(A^{1/2}),
\label{est:Ru-3}
\end{equation}
and
\begin{equation}
\langle Ru,(B^{2}-\lambda A)u\rangle\geq\min\left\{\mu_{2}^{2}\mu_{1}\mu_{3},\lambda_{0}\right\}\cdot|u|^{2}
\qquad
\forall u\in D(B).
\label{th:Ru-below}
\end{equation}

Indeed, since $|\langle u,e_{0}\rangle|\leq|u|$, from (\ref{defn:mu-2}) it follows that
\begin{eqnarray*}
|A^{1/2}Ru| & = & \left|A^{1/2}u-2\langle u,e_{0}\rangle A^{1/2}e_{0}\right| \\[0.5ex]
& \leq & |A^{1/2}u|+2|u|\cdot|A^{1/2}e_{0}|  \\[0.5ex]
& \leq  & |A^{1/2}u|+\frac{2}{\mu_{2}^{1/2}}|A^{1/2}u|\cdot|A^{1/2}e_{0}|,
\end{eqnarray*}
and therefore
\begin{equation}
|\langle Ru,Au\rangle|\leq|A^{1/2}Ru|\cdot|A^{1/2}u|\leq\left(1+\frac{2}{\mu_{2}^{1/2}}|A^{1/2}e_{0}|\right)|A^{1/2}u|^{2}.
\nonumber
\end{equation}

At this point, (\ref{est:Ru-3}) follows from (\ref{defn:delta}).

In order to prove (\ref{th:Ru-below}), we write $u$ as $u_{+}+u_{-}$ according to the decomposition (\ref{defn:H+-}). Since $H_{-}$ and $H_{+}$ are invariant subspaces for $B^{2}-\lambda A$, from (\ref{defn:e0}), (\ref{defn:e0+}) and (\ref{est:B2A}) it follows that
\begin{eqnarray*}
\langle Ru,(B^{2}-\lambda A)u\rangle & = & \langle u_{+}-u_{-},(B^{2}-\lambda A)u_{+}+(B^{2}-\lambda A)u_{-}\rangle \\[1ex]
 & = & \langle u_{+},(B^{2}-\lambda A)u_{+}\rangle-\langle u_{-},(B^{2}-\lambda A)u_{-}\rangle \\[1ex]
 & \geq & \mu_{3}|Bu_{+}|^{2}+\lambda_{0}|u_{-}|^{2} \\[1ex]
 & \geq & \mu_{3}\mu_{2}^{2}\mu_{1}|u_{+}|^{2}+\lambda_{0}|u_{-}|^{2},
\end{eqnarray*}
which implies (\ref{th:Ru-below}).

\paragraph{\textmd{\textit{Differential inequality solved by $F(t)$}}}

We prove that
\begin{equation}
F'(t)\leq -4\gamma_{0}F(t)+|f(t)|^{2}
\qquad
\forall t\geq 0.
\label{diff-ineq:F}
\end{equation}

To this end, we compute the time-derivative of $F(t)$ and we exploit (\ref{eqn:duffing-AB}) and (\ref{defn:Pu}). We obtain that
\begin{eqnarray*}
F'(t) & = & -|u'(t)|^{2}+2\gamma_{0}\langle Pu'(t),u'(t)\rangle \\[1ex]
& & +\langle u'(t),f(t)\rangle+2\gamma_{0}\langle Pu(t),f(t)\rangle \\[1ex]
& & -2\gamma_{0}\left(|Bu(t)|^{2}-\lambda|A^{1/2}u(t)|^{2}+|A^{1/2}u(t)|^{4}\right) \\[1ex]
& & -2\gamma_{0}\delta\langle Ru(t),(B^{2}-\lambda A)u(t)\rangle-2\gamma_{0}\delta|A^{1/2}u(t)|^{2}\cdot\langle Ru(t),Au(t)\rangle.
\end{eqnarray*}

Let $L_{1}$, $L_{2}$, $L_{3}$, $L_{4}$ denote the terms of the four lines. From (\ref{est:pu,u-u,u}) and the first condition in (\ref{defn:gamma0}) we obtain that
\begin{equation}
L_{1}\leq -|u'(t)|^{2}+2\gamma_{0}(1+\delta)|u'(t)|^{2}\leq -\left(\frac{1}{2}+3\gamma_{0}\right)|u'(t)|^{2}.
\nonumber
\end{equation}

From (\ref{est:pu,u-u,u}) and (\ref{est:Pu,Pu-Pu,u}) we obtain that
\begin{eqnarray*}
L_{2} & \leq & \frac{1}{2}|u'(t)|^{2}+\frac{1}{2}|f(t)|^{2}+2\gamma_{0}^{2}|Pu(t)|^{2}+\frac{1}{2}|f(t)|^{2} \\[1ex]
& \leq & \frac{1}{2}|u'(t)|^{2}+2\gamma_{0}^{2}(1+\delta)^{2}|u(t)|^{2}+|f(t)|^{2}.
\end{eqnarray*} 

Finally, from (\ref{th:Ru-below}), (\ref{est:Ru-3}) and the second condition in (\ref{defn:gamma0}) we obtain that
\begin{eqnarray*}
L_{4} & \leq & -2\gamma_{0}\delta\min\left\{\mu_{2}^{2}\mu_{1}\mu_{3},\lambda_{0}\right\}\cdot|u(t)|^{2}+\gamma_{0}|A^{1/2}u(t)|^{4}  \\[0.5ex]
& \leq & -2\gamma_{0}^{2}(5+\delta)(1+\delta)\cdot|u(t)|^{2}+\gamma_{0}|A^{1/2}u(t)|^{4}.
\end{eqnarray*}

Plugging all these estimates into the expression for $F'(t)$, and keeping (\ref{est:F-above}) into account, we deduce (\ref{diff-ineq:F}).

\paragraph{\textmd{\textit{Conclusion}}}

Integrating the differential inequality (\ref{diff-ineq:F}), and letting $t\to+\infty$, we obtain (\ref{th:ultimate}) with $M_{1}:=(4\gamma_{0})^{-1}$. Finally, from (\ref{est:F-below}) and (\ref{est:EAB}) we obtain that
\begin{equation}
|u'(t)|^{2}+|Bu(t)|^{2}\leq\lambda^{2}+4F(t)
\qquad
\forall t\geq 0,
\nonumber
\end{equation}
and therefore (\ref{th:ultimate-u}) with $M_{2}:=\lambda^{2}$ and $M_{3}:=1/\gamma_{0}$ is a consequence of (\ref{th:ultimate}).\qed


\subsection{Proof of Proposition~\ref{prop:unstable}}

\paragraph{\textmd{\textit{Choice of parameters}}}

According to the direct orthogonal sum (\ref{defn:H+-}), any solution $u(t)$ to (\ref{eqn:duffing-AB}) is the sum of a low-frequency component $u_{-}(t)\in H_{-}$ and a high-frequency component $u_{+}(t)\in H_{+}$. Let  $L$ denote as usual the operator $B^{2}-\lambda A$.

The high-frequency component is a solution to
\begin{equation}
u_{+}''(t)+u_{+}'(t)+Lu_{+}(t)+|A^{1/2}u_{+}(t)|^{2}(Au_{+}(t))_{+}=\psi_{1}(t)+\psi_{2}(t),
\label{eqn:u+}
\end{equation}
where 
\begin{eqnarray}
\psi_{1}(t) & := & -|A^{1/2}u_{+}(t)|^{2}(Au_{-}(t))_{+}-2\langle A^{1/2}u_{+}(t),A^{1/2}u_{-}(t)\rangle(Au_{-}(t))_{+} 
\nonumber  \\[1ex]
&& -2\langle A^{1/2}u_{+}(t),A^{1/2}u_{-}(t)\rangle(Au_{+}(t))_{+}-|A^{1/2}u_{-}(t)|^{2}(Au_{+}(t))_{+},
\label{defn:psi-1}
\end{eqnarray}
and
\begin{equation}
\psi_{2}(t):=f_{+}(t)-|A^{1/2}u_{-}(t)|^{2}(Au_{-}(t))_{+}.
\label{defn:psi-2}
\end{equation}

The low-frequency component is a solution to
\begin{equation}
u_{-}''(t)+u_{-}'(t)-\lambda_{0}u_{-}(t)=\psi_{3}(t),
\label{eqn:u-}
\end{equation}
where 
\begin{equation}
\psi_{3}(t):=f_{-}(t)-|A^{1/2}u(t)|^{2}(Au(t))_{-}.
\label{defn:psi-3}
\end{equation}

We recall that $H_{-}$ and $H_{+}$ are not necessarily invariant subspaces for $A$, and for this reason we have to deal with terms of the form $(Au_{\pm}(t))_{\pm}$ in the previous equations.

Now let us set
\begin{equation}
\gamma_{1}:=\min\left\{\frac{1}{24},\frac{\mu_{2}^{2}\mu_{1}\mu_{3}}{14}\right\},
\label{defn:gamma-1}
\end{equation}
let us consider the two constants
\begin{equation}
\Gamma_{1}:=\left(\frac{2}{\gamma_{1}}\frac{1}{\mu_{2}\mu_{1}\mu_{3}}\right)^{3/2},
\qquad\qquad
\Gamma_{2}:=2^{10}\,\Gamma_{1}\cdot|Ae_{0}|^{6},
\nonumber
\end{equation}
and let us choose $\beta_{0}>0$ small enough so that
\begin{equation}
2\cdot 24^{2}\beta_{0}^{4}\cdot|Ae_{0}|^{2}\leq\frac{\gamma_{1}}{2}\mu_{1}\mu_{3},
\label{defn:beta0-1}
\end{equation}
and
\begin{equation}
16\beta_{0}^{2}\cdot|A^{1/2}e_{0}|^{4}+4\Gamma_{2}\beta_{0}^{8}\cdot|A^{1/2}e_{0}|\leq \frac{\lambda_{0}}{2},
\label{defn:beta0-2}
\end{equation}

We claim that, whenever $f(t)$ and $u(t)$ satisfy the assumptions in the left-hand side of (\ref{th:unstable}) with this value of $\beta_{0}$, the solution $u(t)$ satisfies the estimates in the right-hand side of (\ref{th:unstable}) for a suitable constant $M_{4}$ independent of $u(t)$ and $f(t)$. In the sequel we always assume that $f(t)$ and $u(t)$ satisfy the estimates in the left-hand side of (\ref{th:unstable}).

\paragraph{\textmd{\textit{Estimate on right-hand sides}}}

From the assumptions in the left-hand side of (\ref{th:unstable}) we deduce that there exists $t_{0}\geq 0$ such that
\begin{eqnarray}
& |u(t)|\leq 2\beta_{0}
\qquad\qquad
\forall t\geq t_{0}, &
\label{est:b0-u-} \\[1ex]
& |f(t)|\leq 2
\qquad\qquad
\forall t\geq t_{0}. &
\label{est:f}
\end{eqnarray}

We prove that for every $t\geq t_{0}$ it turns out that
\begin{eqnarray}
& |\psi_{1}(t)|\leq 24\beta_{0}^{2}\cdot|Ae_{0}|\cdot |Au_{+}(t)| &
\label{est:psi-1} \\[1ex]
& |\psi_{2}(t)|\leq|f_{+}(t)|+4\beta_{0}^{2}|Ae_{0}|^{2}\cdot|u_{-}(t)|, &
\label{est:psi-2}  \\[1ex]
& |\psi_{3}(t)|\leq|f_{-}(t)|+16\beta_{0}^{2}|A^{1/2}e_{0}|^{4}\cdot|u_{-}(t)|+4|A^{1/2}e_{0}|\cdot|A^{1/2}u_{+}(t)|^{3}. &
\label{est:psi-3}
\end{eqnarray}

To begin with, from (\ref{est:b0-u-}) we deduce that for every $t\geq t_{0}$ it tuns out that
\begin{eqnarray}
& |u_{-}(t)|^{2}\leq 4\beta_{0}^{2}, & 
\label{est:u-1}   \\[1ex]
& |Au_{-}(t)|^{2}\leq 4\beta_{0}^{2}|Ae_{0}|^{2}, & 
\label{est:u-2}   \\[1ex]
& |A^{1/2}u_{-}(t)|^{2}=\langle u_{-}(t),Au_{-}(t)\rangle\leq 4\beta_{0}^{2}|Ae_{0}|. & \label{est:u-3}
\end{eqnarray}

Exploiting (\ref{est:b0-u-}), and (\ref{est:u-1}) through (\ref{est:u-3}), we estimate the four terms in the right-hand side of (\ref{defn:psi-1}), which for the sake of shortness we denote by $T_{1}$, $T_{2}$, $T_{3}$, $T_{4}$. For the first term it turns out that
\begin{equation}
|T_{1}| \leq  |u_{+}(t)|\cdot|Au_{+}(t)|\cdot|Au_{-}(t)| \leq 2\beta_{0}\cdot|Au_{+}(t)|\cdot 2\beta_{0}|Ae_{0}|.
\nonumber
\end{equation}

For the second term it turns out that
\begin{equation}
|T_{2}| \leq 2|Au_{+}(t)|\cdot|u_{-}(t)|\cdot|Au_{-}(t)| \leq 2|Au_{+}(t)|\cdot 2\beta_{0}\cdot 2\beta_{0}|Ae_{0}|.
\nonumber
\end{equation}

For the third term it turns out that
\begin{equation}
|T_{3}| \leq  2|u_{+}(t)|\cdot|Au_{-}(t)|\cdot|(Au_{+}(t))_{+}| \leq 4\beta_{0}\cdot 2\beta_{0}|Ae_{0}|\cdot |Au_{+}(t)|.
\nonumber
\end{equation}

For the fourth term it turns out that
\begin{equation}
|T_{4}| \leq  |A^{1/2}u_{-}(t)|^{2}\cdot|(Au_{+}(t))_{+}| \leq 4\beta_{0}^{2}|Ae_{0}|\cdot|Au_{+}(t)|.
\nonumber
\end{equation}

Plugging the last four inequalities into (\ref{defn:psi-1}) we deduce (\ref{est:psi-1}).

Let us consider now $\psi_{2}(t)$. From (\ref{est:u-3}) we obtain that
\begin{eqnarray*}
\left||A^{1/2}u_{-}(t)|^{2}(Au_{-}(t))_{+}\right| & \leq & |A^{1/2}u_{-}(t)|^{2}\cdot|Au_{-}(t)| \\[0.5ex] 
& \leq & 4\beta_{0}^{2}|Ae_{0}|\cdot|u_{-}(t)|\cdot|Ae_{0}|.
\end{eqnarray*}

Plugging this estimate into (\ref{defn:psi-2}) we deduce (\ref{est:psi-2}).

As for $\psi_{3}(t)$, we observe that
\begin{equation}
(Au(t))_{-}=\langle Au(t),e_{0}\rangle e_{0}=\langle A^{1/2}u(t),A^{1/2}e_{0}\rangle e_{0},
\nonumber
\end{equation}
and therefore
\begin{eqnarray*}
\left||A^{1/2}u(t)|^{2}(Au(t))_{-}\right| & \leq & |A^{1/2}e_{0}|\cdot|A^{1/2}u(t)|^{3} \\
& \leq & |A^{1/2}e_{0}|\cdot\left(|A^{1/2}u_{-}(t)|+|A^{1/2}u_{+}(t)|\right)^{3}  \\[1ex]
& \leq & |A^{1/2}e_{0}|\cdot\left(4|u_{-}(t)|^{3}\cdot|A^{1/2}e_{0}|^{3}+4|A^{1/2}u_{+}(t)|^{3}\strut\right)  \\[1ex]
& \leq & 4|A^{1/2}e_{0}|^{4}\cdot 4\beta_{0}^{2}|u_{-}(t)|+4|A^{1/2}e_{0}|\cdot|A^{1/2}u_{+}(t)|^{3}.
\end{eqnarray*}

Plugging this estimate into (\ref{defn:psi-3}) we deduce (\ref{est:psi-3}).

\paragraph{\textmd{\textit{Estimates on the high-frequency component}}}

We prove that
\begin{equation}
\limsup_{t\to +\infty}\left(|u_{+}'(t)|^{2}+|L^{1/2}u_{+}(t)|^{2}\right)\leq\frac{2}{\gamma_{1}}\limsup_{t\to +\infty}|\psi_{2}(t)|^{2},
\label{th:est-u+}
\end{equation}
and
\begin{equation}
\limsup_{t\to +\infty}|A^{1/2}u_{+}(t)|^{3}\leq 16\Gamma_{1}\cdot\limsup_{t\to +\infty}|f(t)|+\Gamma_{2}\beta_{0}^{8}\cdot\limsup_{t\to +\infty}|u_{-}(t)|.
\label{th:est-Lu+}
\end{equation}

To this end, we consider the energy
\begin{eqnarray*}
F_{+}(t) & := & \frac{1}{2}|u_{+}'(t)|^{2}+\frac{1}{2}|L^{1/2}u_{+}(t)|^{2}+\frac{1}{4}|A^{1/2}u_{+}(t)|^{4} \\[1ex]
& & \mbox{}+2\gamma_{1}\langle u_{+}(t),u_{+}'(t)\rangle+\gamma_{1}|u_{+}(t)|^{2}.
\end{eqnarray*}

Since
\begin{equation}
\left|2\langle u_{+}(t),u_{+}'(t)\rangle\right|\leq|u_{+}(t)|^{2}+|u_{+}'(t)|^{2},
\nonumber
\end{equation}
this energy can be estimated from below by
\begin{equation}
F_{+}(t)\geq\left(\frac{1}{2}-\gamma_{1}\right)|u_{+}'(t)|^{2}+\frac{1}{2}|L^{1/2}u_{+}(t)|^{2},
\label{est:F+below}
\end{equation}
and from above by
\begin{equation}
F_{+}(t)\leq\left(\frac{1}{2}+\gamma_{1}\right)|u_{+}'(t)|^{2}+\frac{1}{2}|L^{1/2}u_{+}(t)|^{2}+\frac{1}{4}|A^{1/2}u_{+}(t)|^{4}+2\gamma_{1}|u_{+}(t)|^{2}.
\label{est:F+above}
\end{equation}

Let us compute the time-derivative of $F_{+}(t)$. Keeping (\ref{eqn:u+}) into account, we obtain that (for the sake of shortness, we omit here the explicit dependence on $t$)
\begin{eqnarray}
F_{+}' & = & -(1-2\gamma_{1})|u_{+}'|^{2}-2\gamma_{1}|L^{1/2}u_{+}|^{2}-2\gamma_{1}|A^{1/2}u_{+}|^{4} \nonumber \\[1ex]
& & +\langle u_{+}',\psi_{1}\rangle+\langle u_{+}',\psi_{2}\rangle+2\gamma_{1}\langle u_{+},\psi_{1}\rangle+2\gamma_{1}\langle u_{+},\psi_{2}\rangle.
\label{deriv:F+}
\end{eqnarray}

We point out that in the computation we exploited identities such as 
\begin{equation}
\langle u_{+}'(t),(Au_{+}(t))_{+}\rangle=\langle u_{+}'(t),Au_{+}(t)\rangle,
\nonumber
\end{equation}
and 
\begin{equation}
\langle u_{+}(t),(Au_{+}(t))_{+}\rangle=\langle u_{+}(t),Au_{+}(t)\rangle=|A^{1/2}u_{+}(t)|^{2}.
\nonumber
\end{equation}

Now we estimate the terms in (\ref{deriv:F+}). As for the terms with $\psi_{2}(t)$, we simply observe that
\begin{equation}
\langle u_{+}'(t),\psi_{2}(t)\rangle\leq\frac{1}{2}|u_{+}'(t)|^{2}+\frac{1}{2}|\psi_{2}(t)|^{2},
\nonumber
\end{equation}
and
\begin{equation}
2\gamma_{1}\langle u_{+}(t),\psi_{2}(t)\rangle\leq 2\gamma_{1}^{2}|u_{+}(t)|^{2}+\frac{1}{2}|\psi_{2}(t)|^{2}.
\nonumber
\end{equation}

In a similar way, for the terms with $\psi_{1}(t)$ we observe that
\begin{equation}
\langle u_{+}'(t),\psi_{1}(t)\rangle\leq\frac{1}{4}|u_{+}'(t)|^{2}+|\psi_{1}(t)|^{2},
\nonumber
\end{equation}
and
\begin{equation}
2\gamma_{1}\langle u_{+}(t),\psi_{1}(t)\rangle\leq \gamma_{1}^{2}|u_{+}(t)|^{2}+|\psi_{1}(t)|^{2}.
\nonumber
\end{equation}

Plugging all these estimates into (\ref{deriv:F+}), and keeping (\ref{est:psi-1}) into account, we conclude that
\begin{eqnarray*}
F_{+}'(t) & \leq & -\left(\frac{1}{4}-2\gamma_{1}\right)|u_{+}'(t)|^{2}-2\gamma_{1}|L^{1/2}u_{+}(t)|^{2}-2\gamma_{1}|A^{1/2}u_{+}(t)|^{4} \nonumber \\[1ex]
& & +3\gamma_{1}^{2}|u_{+}(t)|^{2}+2\cdot 24^{2}\beta_{0}^{4}\cdot|Ae_{0}|^{2}\cdot|Au_{+}(t)|^{2}+|\psi_{2}(t)|^{2}
\end{eqnarray*}
for every $t\geq t_{0}$. Now we claim that
\begin{equation}
F_{+}'(t)\leq -2\gamma_{1}F_{+}(t)+|\psi_{2}(t)|^{2}
\qquad\qquad
\forall t\geq t_{0}.
\label{F+:diff-ineq}
\end{equation}

To this end, keeping (\ref{est:F+above}) into account, it is enough to show that
\begin{eqnarray}
\left(\frac{1}{4}-3\gamma_{1}-2\gamma_{1}^{2}\right)|u_{+}'(t)|^{2}+\frac{3\gamma_{1}}{2}|A^{1/2}u_{+}(t)|^{4}+\gamma_{1}|L^{1/2}u_{+}(t)|^{2} & & 
\nonumber  \\[1ex]
-7\gamma_{1}^{2}|u_{+}(t)|^{2}-2\cdot 24^{2}\beta_{0}^{4}\cdot|Ae_{0}|^{2}\cdot|Au_{+}(t)|^{2} & \geq & 0.
\label{F+:diff-ineq+}
\end{eqnarray}

The coefficient of $|u_{+}'(t)|^{2}$ is positive because $\gamma_{1}\leq 1/24$. Moreover, from (\ref{est:L-B-A}), (\ref{est:B2A}), and the second condition in (\ref{defn:gamma-1}) it turns out that
\begin{equation}
\frac{\gamma_{1}}{2}|L^{1/2}u_{+}(t)|^{2}\geq\frac{\gamma_{1}}{2}\mu_{2}^{2}\mu_{1}\mu_{3}\cdot|u_{+}(t)|^{2}\geq 7\gamma_{1}^{2}|u_{+}(t)|^{2},
\nonumber
\end{equation}
while from (\ref{defn:beta0-1}) it follows that
\begin{equation}
\frac{\gamma_{1}}{2}|L^{1/2}u_{+}(t)|^{2}\geq\frac{\gamma_{1}}{2}\mu_{1}\mu_{3}\cdot|Au_{+}(t)|^{2}\geq 2\cdot 24^{2}\beta_{0}^{4}\cdot|Ae_{0}|^{2}\cdot|Au_{+}(t)|^{2}.
\nonumber
\end{equation}

This completes the proof of (\ref{F+:diff-ineq+}), and therefore also of (\ref{F+:diff-ineq}). Integrating this differential inequality we deduce that
\begin{equation}
\limsup_{t\to +\infty}F_{+}(t)\leq\frac{1}{2\gamma_{1}}\limsup_{t\to +\infty}|\psi_{2}(t)|^{2}
\nonumber
\end{equation}

Since $\gamma_{1}\leq 1/4$, this inequality and (\ref{est:F+below}) imply (\ref{th:est-u+}). 

It remains to prove (\ref{th:est-Lu+}). As usual, from (\ref{est:L-B-A}) and (\ref{est:B2A}) we obtain that
\begin{equation}
|L^{1/2}u_{+}(t)|^{2}\geq\mu_{3}\mu_{1}\mu_{2}|A^{1/2}u_{+}(t)|^{2},
\nonumber
\end{equation}
and therefore from (\ref{th:est-u+}) we deduce that
\begin{eqnarray}
\limsup_{t\to +\infty}|A^{1/2}u_{+}(t)|^{3} & \leq & \left(\frac{1}{\mu_{2}\mu_{1}\mu_{3}}\right)^{3/2}\limsup_{t\to +\infty}|L^{1/2}u_{+}(t)|^{3}  
\nonumber   \\[1ex]
& \leq & \Gamma_{1}\limsup_{t\to +\infty}|\psi_{2}(t)|^{3}.
\label{est:A+3}
\end{eqnarray}

On the other hand, from (\ref{est:psi-2}), (\ref{est:u-1}) and (\ref{est:f}) we obtain that
\begin{eqnarray*}
|\psi_{2}(t)|^{3} & \leq & 4|f(t)|^{3}+4\cdot 64\beta_{0}^{6}\cdot|Ae_{0}|^{6}\cdot|u_{-}(t)|^{3} \\[1ex]
& \leq & 16|f(t)|+2^{10}\beta_{0}^{8}\cdot|Ae_{0}|^{6}\cdot|u_{-}(t)|
\end{eqnarray*}
for every $t\geq t_{0}$. Plugging this estimate into (\ref{est:A+3}), and letting $t\to +\infty$, we obtain (\ref{th:est-Lu+}).
 
\paragraph{\textmd{\textit{Estimates on the low-frequency component}}}

We prove that
\begin{equation}
\limsup_{t\to +\infty}|u_{-}(t)|\leq\frac{2\left(1+64\,\Gamma_{1}|A^{1/2}e_{0}|\right)}{\lambda_{0}}\limsup_{t\to +\infty}|f(t)|.
\label{th:est-u-}
\end{equation}

To this end, we recall that $u_{-}(t)$ is a solution to (\ref{eqn:u-}). Since both the solution $u_{-}(t)$ and the right-hand side $\psi_{3}(t)$ are bounded for positive times, this equation fits in the framework of Lemma~\ref{lemma:ODE-limsup}, from which we deduce that
\begin{equation}
\limsup_{t\to +\infty}|u_{-}(t)|\leq\frac{1}{\lambda_{0}}\limsup_{t\to +\infty}|\psi_{3}(t)|.
\label{est:u-psi-3}
\end{equation}

On the other hand, from (\ref{est:psi-3}) and (\ref{th:est-Lu+}) we know that
\begin{eqnarray}
\limsup_{t\to +\infty}|\psi_{3}(t)| & \leq & \limsup_{t\to +\infty}|f(t)|+16\beta_{0}^{2}|A^{1/2}e_{0}|^{4}\cdot\limsup_{t\to +\infty}|u_{-}(t)| 
\nonumber  \\
& & +4|A^{1/2}e_{0}|\cdot\limsup_{t\to +\infty}|A^{1/2}u_{+}(t)|^{3} 
\nonumber  \\
& \leq & \left(1+64\,\Gamma_{1}|A^{1/2}e_{0}|\right)\limsup_{t\to +\infty}|f(t)|  
\nonumber   \\
& & +\left(16\beta_{0}^{2}\cdot|A^{1/2}e_{0}|^{4}+4\Gamma_{2}\beta_{0}^{8}\cdot|A^{1/2}e_{0}|\right)\limsup_{t\to +\infty}|u_{-}(t)|.
\label{limsup:psi-3}
\end{eqnarray}

Plugging this inequality into (\ref{est:u-psi-3}), and keeping the smallness assumption (\ref{defn:beta0-2}) into account, we conclude that
\begin{equation}
\limsup_{t\to +\infty}|u_{-}(t)|\leq\frac{1+64\,\Gamma_{1}|A^{1/2}e_{0}|}{\lambda_{0}}\limsup_{t\to +\infty}|f(t)|+\frac{1}{2}\limsup_{t\to +\infty}|u_{-}(t)|,
\nonumber
\end{equation}
which implies (\ref{th:est-u-}).

\paragraph{\textmd{\textit{Conclusion}}}

We prove that $u(t)$ satisfies the estimate in the right-hand side of (\ref{th:unstable}) for a suitable constant $M_{4}$. In the sequel $c_{1}$, \ldots, $c_{5}$ denote suitable constants, all independent of the solution and of the forcing term.

From (\ref{th:est-u-}) we know that
\begin{equation}
\limsup_{t\to +\infty}|u_{-}(t)|\leq c_{1}\limsup_{t\to +\infty}|f(t)|.
\label{est:u-psi-3-post}
\end{equation}

Plugging this estimate into (\ref{limsup:psi-3}) we deduce that
\begin{equation}
\limsup_{t\to +\infty}|\psi_{3}(t)|\leq c_{2}\limsup_{t\to +\infty}|f(t)|.
\nonumber
\end{equation}

Applying again Lemma~\ref{lemma:ODE-limsup} to equation (\ref{eqn:u-}) we deduce that
\begin{equation}
\limsup_{t\to +\infty}|u_{-}'(t)|\leq 2\limsup_{t\to +\infty}|\psi_{3}(t)|\leq c_{3}\limsup_{t\to +\infty}|f(t)|.
\nonumber
\end{equation}

Plugging (\ref{est:u-psi-3-post}) into (\ref{est:psi-2}) we obtain that
\begin{equation}
\limsup_{t\to +\infty}|\psi_{2}(t)|\leq c_{4}\limsup_{t\to +\infty}|f(t)|,
\nonumber
\end{equation}
and therefore from (\ref{th:est-u+}) we deduce that
\begin{equation}
\limsup_{t\to +\infty}\left(|u_{+}'(t)|^{2}+|L^{1/2}u_{+}(t)|^{2}\right)\leq c_{5}\limsup_{t\to +\infty}|f(t)|^{2}.
\nonumber
\end{equation}

Recalling (\ref{est:L-B-A}), all these estimates imply the conclusion.\qed


\subsection{Proof of Proposition~\ref{prop:stable}}

\paragraph{\textmd{\textit{Choice of parameters}}}

We can assume, up to replacing $\eta$ with a smaller positive real number, that
\begin{equation}
\eta<\frac{1}{4}\sigma_{0}^{4}{}.
\label{defn:beta-wlog}
\end{equation}

Let us consider the inequality
\begin{equation}
\frac{1}{4}x^{4}-\frac{1}{2}\sigma_{0}^{2}x^{2}< -\frac{\eta}{4}.
\label{defn:x1}
\end{equation}

Due to (\ref{defn:beta-wlog}), the number in the right-hand side is negative but larger than the minimum of the function in the left-hand side. Therefore, the set of solutions to this inequality is the union of two disjoint intervals of the form $(x_{1},x_{2})$ and $(-x_{2},-x_{1})$ for suitable real numbers $0<x_{1}<x_{2}<\sigma_{0}\sqrt{2}$. 

Now let us choose $\gamma_{2}>0$ such that
\begin{equation}
\gamma_{2}\leq\frac{1}{8},
\qquad\qquad
\frac{2\gamma_{2}(1+2\gamma_{2})}{\mu_{2}}\leq(2-\gamma_{2})(\lambda_{2}-\lambda),
\label{defn:gamma-2}
\end{equation}
\begin{equation}
\gamma_{2}\left(\frac{\lambda^{2}}{2}+\frac{2\lambda^{2}}{\mu_{2}^{2}\mu_{1}}+4\sigma_{0}^{2}\right)\leq\frac{\eta}{2},
\label{defn:gamma-2-bis}
\end{equation}
and
\begin{equation}
\frac{\gamma_{2}}{2}\sigma_{0}^{2}+\frac{2\gamma_{2}(1+2\gamma_{2})}{\mu_{2}{}}\leq\sigma_{0}(2-\gamma_{2})x_{1}.
\label{defn:gamma-2-ter}
\end{equation}

Finally, let us choose $\ep_{1}>0$ such that
\begin{equation}
\frac{\ep_{1}^{2}}{2\gamma_{2}^{2}}<\frac{\eta}{4}.
\label{defn:ep-1}
\end{equation}

\paragraph{\textmd{\textit{Estimate at time $T_{0}$}}}

Let us write $u(t)$ in the form $\alpha(t)e_{1}+w(t)$ according to the direct sum (\ref{defn:W+-}). We prove that
\begin{equation}
|\alpha(T_{0})|>x_{1}.
\label{th:at0}
\end{equation}

Indeed, from (\ref{est:E-below-aw}) we obtain that
\begin{equation}
E(t)\geq\frac{1}{2}|u'(t)|^{2}+\frac{1}{4} \alpha^{4}(t)-\frac{1}{2}\sigma_{0}^{2} \alpha^{2}(t)+\frac{1}{2}\frac{\lambda_{2}-\lambda}{\lambda_{2}}|Bw(t)|^{2}.
\nonumber
\end{equation}

Setting $t=T_{0}$, from the assumption that $E(T_{0})<-\eta$ we conclude that
\begin{equation}
\frac{1}{4} \alpha^{4}(T_{0})-\frac{1}{2}\sigma_{0}^{2} \alpha^{2}(T_{0})<-\eta.
\nonumber
\end{equation}

Comparing with (\ref{defn:x1}) we deduce (\ref{th:at0}).

\paragraph{\textmd{\textit{Modified energy and basic estimates from above and below}}}

Due to (\ref{th:at0}) and the symmetry of the problem, in the sequel we can assume, without loss of generality, that $\alpha(T_{0})>x_{1}$. In this case we claim that the solution is eventually close to the stationary point $\sigma_{0}e_{1}$, and for this reason we introduce the modified energy
\begin{equation}
S(t):=E(t)+\frac{1}{4}\sigma_{0}^{4}{}+2\gamma_{2}\langle u(t)-\sigma_{0}e_{1},u'(t)\rangle+\gamma_{2}|u(t)-\sigma_{0}e_{1}|^{2}.
\nonumber
\end{equation}

From the inequality
\begin{equation}
2\langle u(t)-\sigma_{0}e_{1},u'(t)\rangle\leq|u(t)-\sigma_{0}e_{1}|^{2}+|u'(t)|^{2}
\label{cs:gamma-2}
\end{equation}
we deduce that
\begin{eqnarray}
S(t) & \leq & \left(\frac{1}{2}+\gamma_{2}\right)|u'(t)|^{2}+\frac{1}{2}|Bu(t)|^{2}-\frac{\lambda}{2}|A^{1/2}u(t)|^{2}+\frac{1}{4}|A^{1/2}u(t)|^{4} 
\nonumber  \\
& & +\frac{1}{4}\sigma_{0}^{4}{}+2\gamma_{2}|u(t)-\sigma_{0}e_{1}|^{2},
\label{est:S-above}
\end{eqnarray}
and
\begin{equation}
S(t)\geq\left(\frac{1}{2}-\gamma_{2}\right)|u'(t)|^{2}+\frac{1}{2}|Bu(t)|^{2}-\frac{\lambda}{2}|A^{1/2}u(t)|^{2}+\frac{1}{4}|A^{1/2}u(t)|^{4}+\frac{1}{4}\sigma_{0}^{4}{}.
\nonumber
\end{equation}

If we write $u(t)$ in the usual form $\alpha(t)e_{1}+w(t)$, and we keep (\ref{eqn:Bw}) into account, the estimate from below implies that
\begin{equation}
S(t)\geq\left(\frac{1}{2}-\gamma_{2}\right)|u'(t)|^{2}+\frac{1}{4}\left( \alpha^{2}(t)-\sigma_{0}^{2}\right)^{2}{}+\frac{1}{2}\frac{\lambda_{2}-\lambda}{\lambda_{2}}|Bw(t)|^{2}.
\label{est:S-below}
\end{equation}

\paragraph{\textmd{\textit{Modified energy at time $T_{0}$}}}

We prove that
\begin{equation}
S(T_{0})<\frac{1}{4}\sigma_{0}^{4}{}-\frac{\eta}{2}.
\label{th:St0}
\end{equation}

Indeed, the energies $S(t)$ and $E(t)$ satisfy
\begin{eqnarray}
S(t) & = & E(t)+\frac{1}{4}\sigma_{0}^{4}{} 
\nonumber  \\[0.5ex]
& & +2\gamma_{2}\langle u(t)-\sigma_{0}e_{1},u'(t)\rangle+\gamma_{2}|u(t)-\sigma_{0}e_{1}|^{2}.
\label{eqn:S(t)}
\end{eqnarray}

Let $\Lambda$ denote the sum of the two terms in the last line. From (\ref{cs:gamma-2}) we know that
\begin{equation}
\Lambda \leq \gamma_{2}|u'(t)|^{2}+2\gamma_{2}|u(t)-\sigma_{0}e_{1}|^{2}\leq\gamma_{2}\left(|u'(t)|^{2}+4|u(t)|^{2}+4\sigma_{0}^{2}\right).
\label{est:L}
\end{equation}

Setting $t=T_{0}$, from (\ref{est:E-below}) we obtain that
\begin{equation}
|u'(T_{0})|^{2}\leq 2E(T_{0})+\frac{\lambda^{2}}{2}\leq\frac{\lambda^{2}}{2},
\nonumber
\end{equation}
and similarly from (\ref{est:E-below}) and (\ref{est:B2A}) we obtain that
\begin{equation}
|u(T_{0})|^{2}\leq\frac{1}{\mu_{2}^{2}\mu_{1}}|Bu(T_{0})|^{2}\leq\frac{1}{\mu_{2}^{2}\mu_{1}}\left(2E(T_{0})+\frac{\lambda^{2}}{2}\right)\leq\frac{1}{\mu_{2}^{2}\mu_{1}}\frac{\lambda^{2}}{2}.
\label{est:L3-B}
\end{equation}

Replacing the last two inequalities into (\ref{est:L}) we deduce that
\begin{equation}
\Lambda\leq\gamma_{2}\left(\frac{\lambda^{2}}{2}+\frac{2\lambda^{2}}{\mu_{2}^{2}\mu_{1}}+4\sigma_{0}^{2}\right).
\nonumber
\end{equation}

Plugging this estimate into (\ref{eqn:S(t)}), and keeping the smallness assumption (\ref{defn:gamma-2-bis}) into account, we deduce (\ref{th:St0}).

\paragraph{\textmd{\textit{Modified energy and potential well}}}

We show that, for every $t\geq 0$, the following implication holds true:
\begin{equation}
S(t)<\frac{1}{4}\sigma_{0}^{4}{}-\frac{\eta}{4}
\qquad\Longrightarrow\qquad
|\alpha(t)|>x_{1}.
\label{th:S-pot-well}
\end{equation}

Indeed, from (\ref{est:S-below}) we know that
\begin{equation}
S(t)\geq \frac{1}{4}\left( \alpha^{2}(t)-\sigma_{0}^{2}\right)^{2}{},
\nonumber
\end{equation}
and therefore the inequality in the left-hand side of (\ref{th:S-pot-well}) implies that
\begin{equation}
\frac{1}{4} \alpha^{4}(t)-\frac{1}{2}\sigma_{0}^{2} \alpha^{2}(t)<-\frac{\eta}{4}.
\nonumber
\end{equation}

Comparing with (\ref{defn:x1}), this implies that $|\alpha(t)|>x_{1}$.

\paragraph{\textmd{\textit{Differential inequality in the potential-well}}}

We prove that, for every $t\geq 0$, the following implication holds true:
\begin{equation}
\alpha(t)\geq x_{1}
\qquad\Longrightarrow\qquad
S'(t)\leq -2\gamma_{2}^{2}S(t)+|f(t)|^{2}.
\label{S'-pot-well}
\end{equation}

To begin with, we compute the time-derivative of $S(t)$, which turns out to be 
\begin{eqnarray*}
S'(t) & = & -(1-2\gamma_{2})|u'(t)|^{2}+\langle u'(t),f(t)\rangle+2\gamma_{2}\langle u(t)-\sigma_{0}e_{1},f(t)\rangle \\[1ex]
& & -2\gamma_{2}\langle u(t)-\sigma_{0}e_{1},B^{2}u(t)-\lambda Au(t)+|A^{1/2}u(t)|^{2}Au(t)\rangle
\end{eqnarray*}

From the usual inequalities
\begin{equation}
\langle u'(t),f(t)\rangle\leq\frac{1}{2}|u'(t)|^{2}+\frac{1}{2}|f(t)|^{2}
\nonumber
\end{equation}
and
\begin{equation}
2\gamma_{2}\langle u(t)-\sigma_{0}e_{1},f(t)\rangle\leq 2\gamma_{2}^{2}|u(t)-\sigma_{0}e_{1}|^{2}+\frac{1}{2}|f(t)|^{2}
\nonumber
\end{equation}
we deduce that
\begin{eqnarray*}
S'(t) & \leq & -\left(\frac{1}{2}-2\gamma_{2}\right)|u'(t)|^{2}+|f(t)|^{2}+2\gamma_{2}^{2}|u(t)-\sigma_{0}e_{1}|^{2} 
\nonumber  \\
& & -2\gamma_{2}\langle u(t),B^{2}u(t)-\lambda Au(t)+|A^{1/2}u(t)|^{2}Au(t)\rangle 
\nonumber  \\[0.5ex] 
& & +2\gamma_{2}\sigma_{0}\langle e_{1},B^{2}u(t)-\lambda Au(t)+|A^{1/2}u(t)|^{2}Au(t)\rangle.
\end{eqnarray*}

Keeping (\ref{est:S-above}) into account, inequality (\ref{S'-pot-well}) is proved if we show that
\begin{eqnarray}
\left(\frac{1}{2}-2\gamma_{2}-\gamma_{2}^{2}-2\gamma_{2}^{3}\right)|u'(t)|^{2}-2\gamma_{2}^{2}(1+2\gamma_{2})|u(t)-\sigma_{0}e_{1}|^{2} & & 
\nonumber  \\
+\gamma_{2}(2-\gamma_{2})\left(|Bu(t)|^{2}-\lambda|A^{1/2}u(t)|^{2}\right)+\gamma_{2}\left(2-\frac{\gamma_{2}}{2}\right)|A^{1/2}u(t)|^{4} & & 
\nonumber  \\
-2\gamma_{2}\sigma_{0}\langle e_{1},B^{2}u(t)-\lambda Au(t)+|A^{1/2}u(t)|^{2}Au(t)\rangle-\frac{\gamma_{2}^{2}}{2}\sigma_{0}^{4}{} & \geq & 0.
\label{est:S'}
\end{eqnarray} 

Now we write $u(t)$ in the usual form $\alpha(t)e_{1}+w(t)$ according to the direct sum (\ref{defn:W+-}), and we estimate the terms in the left-hand side. 
\begin{itemize}

\item  The coefficient of $|u'(t)|^{2}$ is nonnegative because $\gamma_{2}\leq 1/8$.

\item  For the second term we exploit (\ref{defn:mu-2}) and (\ref{A=a+w}), obtaining that
\begin{eqnarray*}
|u(t)-\sigma_{0}e_{1}|^{2} & \leq & \frac{1}{\mu_{2}}|A^{1/2}(u(t)-\sigma_{0}e_{1})|^{2} \\
& = & \frac{1}{\mu_{2}}\left\{(\alpha(t)-\sigma_{0})^{2}{}+|A^{1/2}w(t)|^{2}\right\}.
\end{eqnarray*}

\item  For the third term we exploit (\ref{eqn:Bu-Au}), (\ref{est:Bw}) and (\ref{hp:e1-perp}), and we obtain that
\begin{eqnarray*}
|Bu(t)|^{2}-\lambda|A^{1/2}u(t)|^{2} & = & -\sigma_{0}^{2} \alpha^{2}(t){}+|Bw(t)|^{2}-\lambda|A^{1/2}w(t)|^{2}  \\[1ex]
& \geq & -\sigma_{0}^{2} \alpha^{2}(t){}+(\lambda_{2}-\lambda)|A^{1/2}w(t)|^{2}.
\end{eqnarray*}

\item  We expand the fourth term according to (\ref{A=a+w}). 

\item  Finally, from (\ref{eqn:w-e1}), (\ref{hp:e1}) and (\ref{defn:sigma-0}) we know that
\begin{equation}
\langle e_{1},B^{2}u(t)-\lambda Au(t)\rangle=\alpha(t)\lambda_{1}{}-\alpha(t)\lambda{}=-\sigma_{0}^{2}\alpha(t){},
\nonumber
\end{equation}
while from (\ref{A=a+w}) we know that
\begin{equation}
|A^{1/2}u(t)|^{2}\cdot\langle e_{1},Au(t)\rangle=\left( \alpha^{2}(t){}+|A^{1/2}w(t)|^{2}\right)\cdot\alpha(t){},
\nonumber
\end{equation}
and therefore the scalar product in the fifth term is equal to
\begin{equation}
-\sigma_{0}^{2}\alpha(t){}+ \alpha^{3}(t){}+\alpha(t){}\cdot|A^{1/2}w(t)|^{2}.
\nonumber
\end{equation}

\end{itemize}

Keeping all these equalities and inequalities into account, we obtain that (\ref{est:S'}) holds true if we show that
\begin{equation}
k_{1}(t)|A^{1/2}w(t)|^{2}+k_{2}(t)|A^{1/2}w(t)|^{4}+2k_{3}(t){}|A^{1/2}w(t)|^{2}+k_{4}(t){}\geq 0,
\label{est:S'-final}
\end{equation}
where
\begin{equation}
k_{1}(t):=(2-\gamma_{2})(\lambda_{2}-\lambda)-\frac{2\gamma_{2}(1+2\gamma_{2})}{\mu_{2}},
\nonumber
\end{equation}
\begin{equation}
k_{2}(t):=2-\frac{\gamma_{2}}{2},   
\qquad\qquad
k_{3}(t):=\alpha(t)\left\{\left(2-\frac{\gamma_{2}}{2}\right)\alpha(t)-\sigma_{0}\right\},
\nonumber
\end{equation}
and
\begin{eqnarray*}
k_{4}(t) & := & -\frac{2\gamma_{2}(1+2\gamma_{2})}{\mu_{2}{}}\left(\alpha(t)-\sigma_{0}\right)^{2}-(2-\gamma_{2})\sigma_{0}^{2} \alpha^{2}(t)+\left(2-\frac{\gamma_{2}}{2}\right) \alpha^{4}(t)  \\[0.5ex]
& & +2\sigma_{0}^{3}\alpha(t)-2\sigma_{0} \alpha^{3}(t)-\frac{\gamma_{2}}{2}\sigma_{0}^{4}  \\
& = & \left\{\left(2-\frac{\gamma_{2}}{2}\right) \alpha^{2}(t)+\sigma_{0}(2-\gamma_{2})\alpha(t)-\frac{\gamma_{2}}{2}\sigma_{0}^{2}-\frac{2\gamma_{2}(1+2\gamma_{2})}{\mu_{2}{}}\right\}\left(\alpha(t)-\sigma_{0}\right)^{2}.
\end{eqnarray*}

Now we exploit the assumption that $\alpha(t)\geq x_{1}$. When this is the case, from the smallness assumptions (\ref{defn:gamma-2}) it follows that 
\begin{equation}
k_{1}(t)\geq 0,
\qquad
k_{2}(t)\geq 1,
\qquad
k_{3}(t)\geq\alpha(t)(\alpha(t)-\sigma_{0}),
\nonumber
\end{equation}
while from the smallness assumption (\ref{defn:gamma-2-ter}) it follows that
\begin{equation}
k_{4}(t)\geq \alpha^{2}(t)(\alpha(t)-\sigma_{0})^{2}.
\nonumber
\end{equation}

As a consequence, the left-hand side of (\ref{est:S'-final}) is greater than or equal to
\begin{equation}
\left(|A^{1/2}w(t)|^{2}+\alpha(t)(\alpha(t)-\sigma_{0}){}\right)^{2},
\nonumber
\end{equation}
and therefore it is nonnegative in this regime. This completes the proof of (\ref{S'-pot-well}).

\paragraph{\textmd{\textit{Potential-well argument}}}

We prove that
\begin{equation}
\alpha(t)>x_{1}
\qquad
\forall t\geq T_{0},
\label{th:pot-well}
\end{equation}
which is equivalent to (\ref{th:stable-sign}). To this end, let us set
\begin{equation}
T_{1}:=\sup\left\{t\geq T_{0}: \alpha(\tau)>x_{1}\mbox{ for every }\tau\in[T_{0},t]\right\}.
\nonumber
\end{equation}

We observe that $T_{1}$ is the supremum of an open set containing $t=T_{0}$ because we assumed that $\alpha(T_{0})>0$ after showing (\ref{th:at0}). It follows that $T_{1}$ is well-defined, greater than $T_{0}$, and it satisfies
\begin{equation}
\alpha(t)>x_{1}
\qquad
\forall t\in[T_{0},T_{1}).
\label{hp:T1}
\end{equation}

If $T_{1}=+\infty$, then (\ref{th:pot-well}) is proved. Let us assume by contradiction that $T_{1}<+\infty$. Then the maximality of $T_{1}$ implies that $\alpha(T_{1})=x_{1}$. On the other hand, from (\ref{hp:T1}) it follows that (\ref{S'-pot-well}) holds true for every $t\in[T_{0},T_{1}]$. Integrating this differential inequality, and recalling that $|f(t)|\leq\ep_{1}$ for every $t\geq T_{0}$, we deduce that
\begin{equation}
S(t)\leq S(T_{0})+\frac{\ep_{1}^{2}}{2\gamma_{2}^{2}}
\qquad
\forall t\in[T_{0},T_{1}].
\nonumber
\end{equation}

Keeping (\ref{th:St0}) and (\ref{defn:ep-1}) into account, this implies that
\begin{equation}
S(T_{1})<\frac{1}{4}\sigma_{0}^{4}{}-\frac{\eta}{4},
\nonumber
\end{equation}
which in turn implies that $|\alpha(T_{1})|>x_{1}$ because of (\ref{th:S-pot-well}). This contradicts the fact that $\alpha(T_{1})=x_{1}$, and completes the proof of (\ref{th:pot-well}).

\paragraph{\textmd{\textit{Conclusion}}}

Since we have established (\ref{th:pot-well}), now we know that the differential inequality in (\ref{S'-pot-well}) holds true for every $t\geq T_{0}$. Integrating this differential inequality, and letting $t\to +\infty$, we deduce that
\begin{equation}
\limsup_{t\to +\infty}S(t)\leq\frac{1}{2\gamma_{2}^{2}}\limsup_{t\to +\infty}|f(t)|^{2}.
\label{th:S-ultimate}
\end{equation}

Now we prove that there exists a constant $\Gamma_{3}$ such that
\begin{equation}
|u'(t)|^{2}+|B(u(t)-\sigma_{0}e_{1})|^{2}\leq\Gamma_{3}\,S(t)
\label{th:est-u-S}
\end{equation}
for all solutions $u(t)$ with $\alpha(t)\geq 0$. If we show this claim, then (\ref{th:S-ultimate}) implies the conclusion in the lower box of (\ref{th:stable}).

In order to prove (\ref{th:est-u-S}), it is enough to observe that
\begin{eqnarray}
|B(u(t)-\sigma_{0}e_{1})|^{2} & = & \left(\alpha(t)-\sigma_{0}\right)^{2}|Be_{1}|^{2}+|Bw(t)|^{2}
\nonumber  \\[0.5ex] 
& = & \lambda_{1}\left(\alpha(t)-\sigma_{0}\right)^{2}{}+|Bw(t)|^{2} 
\nonumber  \\[0.5ex]
& \leq & \frac{\left(\alpha(t)+\sigma_{0}\right)^{2}}{\sigma_{0}^{2}}\cdot\lambda_{1}\left(\alpha(t)-\sigma_{0}\right)^{2}{}+|Bw(t)|^{2}  
\nonumber  \\
& = & \frac{\lambda_{1}}{\sigma_{0}^{2}}\left( \alpha^{2}(t)-\sigma_{0}^{2}\right)^{2}{}+|Bw(t)|^{2},
\label{est:u-se}
\end{eqnarray}
where in the inequality we exploited the assumption that $\alpha(t)\geq 0$. At this point it is enough to observe that, due to (\ref{est:S-below}), the energy $S(t)$ controls both $|u'(t)|^{2}$ and the terms in the right-hand side of (\ref{est:u-se}), up to constants.\qed


\subsection{Proof of Proposition~\ref{prop:asymptotic}}

Let $r(t):=u(t)-v(t)$ denote the difference between two solutions $u(t)$ and $v(t)$ to equation (\ref{eqn:duffing-AB}). This difference is a solution to equation
\begin{equation}
r''(t)+r'(t)+B^{2}r(t)-\lambda Ar(t)=g(t),
\label{eqn:r}
\end{equation}
where
\begin{eqnarray}
g(t) & := & -|A^{1/2}u(t)|^{2}Au(t)+|A^{1/2}v(t)|^{2}Av(t) \nonumber \\[0.5ex]
 & = & -|A^{1/2}u(t)|^{2}Ar(t)-\langle u(t)+v(t),Ar(t)\rangle Av(t).
 \label{eqn:g-r}
\end{eqnarray}

Now we consider separately the  unstable case $\sigma=0$ and the stable cases $\sigma=\pm\sigma_{0}$. The constants $c_{6}$, \ldots, $c_{21}$ in the sequel are independent of $u(t)$ and $v(t)$.

\paragraph{\textmd{\textit{Unstable case}}}

From (\ref{eqn:g-r}) we deduce that
\begin{equation}
|g(t)|\leq\left(|A^{1/2}u(t)|^{2}+|u(t)+v(t)|\cdot|Av(t)|\right)|Ar(t)|.
\label{est:g-u}
\end{equation}
On the other hand, from (\ref{hp:asymptotic}) with $\sigma=0$ and (\ref{est:B2A}) we know that
$$\limsup_{t\to+\infty}|A^{1/2}u(t)|^{2}\leq c_{6}r_{0}^{2}
\quad\quad\mbox{and}\quad\quad
\limsup_{t\to+\infty}|u(t)+v(t)|\cdot|Av(t)|\leq c_{7}r_{0}^{2},$$
and therefore from  (\ref{est:g-u}) we obtain that
\begin{equation}
\limsup_{t\to+\infty}|g(t)|\leq c_{8}r_{0}^{2}\cdot\limsup_{t\to+\infty}|Ar(t)|\leq c_{9}r_{0}^{2}\cdot\limsup_{t\to+\infty}|Br(t)|.
\label{est:gAr}
\end{equation}

Let $L$ denote as usual the operator $B^{2}-\lambda A$, let $r_{-}(t)$ and $r_{+}(t)$ denote the components of $r(t)$ with respect to the direct orthogonal sum (\ref{defn:H+-}), and let $g_{-}(t)$ and $g_{+}(t)$ denote the corresponding components of $g(t)$. As already observed, $L$ is the differential of the functional (\ref{defn:funct-AB}) in the origin.

Since $e_{0}$ is an eigenvector of $L$ with eigenvalue $\lambda_{0}$, the low-frequency component $r_{-}(t)$ is a solution to equation
\begin{equation}
r_{-}''(t)+r_{-}'(t)-\lambda_{0}r_{-}(t)=g_{-}(t),
\label{eqn:r-}
\end{equation}
while the high-frequency component $r_{+}(t)$ is a solution to equation
\begin{equation}
r''_{+}(t)+r'_{+}(t)+Lr_{+}(t)=g_{+}(t).
\label{eqn:r+}
\end{equation}

Equation (\ref{eqn:r-}) is a scalar equation that fits in the framework of Lemma~\ref{lemma:ODE-limsup} with
$$y(t):=r_{-}(t),
\quad\quad
m:=\lambda_{0},
\quad\quad
\psi(t):=g_{-}(t).$$

Indeed, $r_{-}(t)$ is bounded because $u(t)$ and $v(t)$ are bounded, and for analogous reasons also $g_{-}(t)$ is bounded. As a consequence, from  Lemma~\ref{lemma:ODE-limsup} we deduce that
\begin{equation}
\limsup_{t\to+\infty}\left(|r_{-}'(t)|^{2}+|Br_{-}(t)|^{2}\right)\leq c_{10}\limsup_{t\to+\infty}|g_{-}(t)|^{2}.
\label{est:r-g-}
\end{equation}

Now from (\ref{est:L-B-A}) and (\ref{est:B2A}) we know that
\begin{equation}
|L^{1/2}x|^{2}\geq\mu_{3}|Bx|^{2}\geq\mu_{2}^{2}\mu_{1}\mu_{3}|x|^{2}
\qquad
\forall x\in H_{+}\cap D(B),
\nonumber
\end{equation}
and therefore equation (\ref{eqn:r+}) fits in the  framework of Lemma~\ref{lemma:PDE-limsup} with
$$X:=H_{+},
\quad\quad
\mathcal{L}:=L,
\quad\quad
y(t):=r_{+}(t),
\quad\quad
m:=\mu_{2}^{2}\mu_{1}\mu_{3},
\quad\quad
\psi(t):=g_{+}(t).$$

As a consequence, from Lemma~\ref{lemma:PDE-limsup} we deduce that
\begin{eqnarray}
\limsup_{t\to+\infty}\left(|r_{+}'(t)|^{2}+|Br_{+}(t)|^{2}\right) & \leq & c_{11}\limsup_{t\to+\infty}\left(|r_{+}'(t)|^{2}+|L^{1/2}r_{+}(t)|^{2}\right) 
\nonumber \\
 & \leq & c_{12}\limsup_{t\to+\infty}|g_{+}(t)|^{2}.
 \label{est:r+g+}
\end{eqnarray}

Since $|Br(t)|^{2}\leq 2|Br_{-}(t)|^{2}+2|Br_{+}(t)|^{2}$ (we recall that $H_{-}$ and $H_{+}$ are not necessarily invariant subspaces for $B$), from (\ref{est:r-g-}), (\ref{est:r+g+}), and (\ref{est:gAr}) we conclude that
\begin{eqnarray*}
\limsup_{t\to+\infty}\left(|r'(t)|^{2}+|Br(t)|^{2}\right) & \leq & 
2\limsup_{t\to+\infty}\left(|r_{-}'(t)|^{2}+|Br_{-}(t)|^{2}\right) \\
 & & \mbox{} +2\limsup_{t\to+\infty}\left(|r_{+}'(t)|^{2}+|Br_{+}(t)|^{2}\right) \\
 & \leq & c_{13}\limsup_{t\to+\infty}|g(t)|^{2} \\
 & \leq & c_{14}r_{0}^{4}\cdot\limsup_{t\to+\infty}|Br(t)|^{2} \\
 & \leq & c_{14}r_{0}^{4}\cdot\limsup_{t\to+\infty}\left(|r'(t)|^{2}+|Br(t)|^{2}\right).
\end{eqnarray*}

If $r_{0}$ is small enough, the coefficient $c_{14}r_{0}^{4}$ is less than~1. It follows that
\begin{equation}
\lim_{t\to+\infty}\left(|r'(t)|^{2}+|Br(t)|^{2}\right)=0,
\label{th:lim-r}
\end{equation}
which in turn is equivalent to (\ref{th:main-asymptotic}).

\paragraph{\textmd{\textit{Stable case}}}

We assume, without loss of generality, that $\sigma=\sigma_{0}$ (the case $\sigma=-\sigma_{0}$ being symmetric). In order to exploit the smallness of $u(t)-\sigma_{0}e_{1}$ and $v(t)-\sigma_{0}e_{1}$, with some algebra we rewrite (\ref{eqn:r}) in the form
\begin{equation}
r''(t)+r'(t)+B^{2}r(t)-\lambda Ar(t)+\sigma_{0}^{2}{}Ar(t)+2\sigma_{0}^{2}\langle r(t),Ae_{1}\rangle Ae_{1}=\hg(t),
\label{eqn:r-stable}
\end{equation}
where
\begin{eqnarray*}
\hg(t) & := & -\left(|A^{1/2}(u(t)-\sigma_{0}e_{1})|^{2}+
2\langle A(u(t)-\sigma_{0}e_{1}),\sigma_{0}e_{1}\rangle\right)Ar(t) \\
 & & -\langle u(t)+v(t),Ar(t)\rangle A(v(t)-\sigma_{0}e_{1}) \\
 & & -\langle u(t)+v(t)-2\sigma_{0}e_{1},Ar(t)\rangle \sigma_{0}Ae_{1}.
\end{eqnarray*}

Due to (\ref{hp:asymptotic}) with $\sigma=\sigma_{0}$, the forcing term $\hg(t)$ satisfies
\begin{equation}
\limsup_{t\to+\infty}|\hg(t)|\leq
\left(c_{15}r_{0}+c_{16}r_{0}^{2}\right)\limsup_{t\to+\infty}|Ar(t)|.
\label{est:g123}
\end{equation}

Now we observe that (\ref{eqn:r-stable}) can be rewritten in the form
\begin{equation}
r''(t)+r'(t)+Cr(t)=\hg(t),
\label{eqn:r-stable-C}
\end{equation}
where $C$ is the linear operator on $H$ (with domain $D(C)=D(B^{2})$) defined by 
\begin{equation}
Cx=B^{2}x-\lambda Ax+\sigma_{0}^{2}{}Ax+2\sigma_{0}^{2}\langle x,Ae_{1}\rangle Ae_{1},
\nonumber
\end{equation}
which coincides with the differential of the functional (\ref{defn:funct-AB}) in $\sigma_{0}e_{1}$.

We claim that there exists a positive constant $c_{17}$ such that
\begin{equation}
|Bx|^{2}\leq c_{17}|C^{1/2}x|^{2}
\qquad
\forall x\in D(B),
\label{est:C2B}
\end{equation}
which implies in particular that 
\begin{equation}
|C^{1/2}x|^{2}\geq c_{18}|x|^{2}
\qquad
\forall x\in D(B)
\nonumber
\end{equation}
for a suitable positive constant $c_{18}$. To this end, we first observe that
\begin{equation}
|C^{1/2}x|^{2}=|Bx|^{2}-\lambda|A^{1/2}x|^{2}+\sigma_{0}^{2}{}\cdot|A^{1/2}x|^{2}+2\sigma_{0}^{2}\langle x,Ae_{1}\rangle^{2}.
\nonumber
\end{equation}

Then we write $x$ in the form $\alpha e_{1}+w$ according to the direct sum (\ref{defn:W+-}), and from (\ref{eqn:Bu-Au}) and (\ref{est:Bw}) we deduce that
\begin{eqnarray*}
|Bx|^{2}-\lambda|A^{1/2}x|^{2} & = & -\sigma_{0}^{2}\alpha^{2}{}+|Bw|^{2}-\lambda|A^{1/2}w|^{2}  \\[0.5ex]
& \geq & -\sigma_{0}^{2}\alpha^{2}{}+\frac{\lambda_{2}-\lambda}{\lambda_{2}}|Bw|^{2},
\end{eqnarray*}
and 
\begin{equation}
\sigma_{0}^{2}{}\cdot|A^{1/2}x|^{2}+2\sigma_{0}^{2}\langle x,Ae_{1}\rangle^{2}=3\sigma_{0}^{2}\alpha^{2}{}+\sigma_{0}^{2}{}\cdot|A^{1/2}w|^{2}.
\nonumber
\end{equation}

From these inequalities it follows that
\begin{equation}
|C^{1/2}u|^{2}\geq\frac{\lambda_{2}-\lambda}{\lambda_{2}}|Bw|^{2}+\sigma_{0}^{2}{}\cdot|A^{1/2}w|^{2}+2\sigma_{0}^{2}\alpha^{2}{}.
\label{est:L1/2}
\end{equation}

On the other hand from (\ref{B=a+w}) we know that
\begin{equation}
|Bx |^{2}=\lambda_{1}\alpha^{2}{}+|Bw|^{2}.
\label{est:Bu-final}
\end{equation}

Comparing (\ref{est:L1/2}) and (\ref{est:Bu-final}) we deduce (\ref{est:C2B}). At this point, equation (\ref{eqn:r-stable-C}) fits in the framework of Lemma~\ref{lemma:PDE-limsup} with
$$X:=H,
\quad\quad
\mathcal{L}:=C,
\qquad
y(t):=r(t),
\quad\quad
m:=c_{18},
\quad\quad
\psi(t):=\hg(t).$$

As a consequence, from (\ref{est:C2B}), Lemma~\ref{lemma:PDE-limsup}, and (\ref{est:g123}) we deduce that
\begin{eqnarray*}
\limsup_{t\to+\infty}\left(|r'(t)|^{2}+|Br(t)|^{2}\right) & \leq & c_{19}\limsup_{t\to+\infty}\left(|r'(t)|^{2}+|C^{1/2}r(t)|^{2}\right) 
\nonumber \\
 & \leq & c_{20}\limsup_{t\to+\infty}|\hg(t)|^{2}.
\nonumber  \\
 & \leq & 
c_{20}\left(c_{15}r_{0}+c_{16}r_{0}^{2}\right)^{2}\cdot\limsup_{t\to+\infty}|Ar(t)|^{2} 
\\
 & \leq & c_{21}(r_{0}^{2}+r_{0}^{4})\cdot\limsup_{t\to+\infty}\left(|r'(t)|^{2}+|Br(t)|^{2}\right).
\end{eqnarray*}

If $r_{0}$ is small enough, we obtain again (\ref{th:lim-r}), which in turn is equivalent to (\ref{th:main-asymptotic}).\qed


\setcounter{equation}{0}
\section{Proof of Theorem~\ref{thm:main}}\label{sec:proof-main}

\paragraph{\textmd{\textit{Choice of parameters}}}

To begin with, we consider the constants $\beta_{0}$ and $M_{4}$ of Proposition~\ref{prop:unstable}. Then we apply Proposition~\ref{prop:stable} with
\begin{equation}
\eta:=\min\left\{\frac{\gamma_{0}}{4}(1-\delta)\beta_{0}^{2},\frac{1}{8}\sigma_{0}^{4}{}\right\},
\label{defn:beta}
\end{equation}
where $\delta$ and $\gamma_{0}$ are chosen in (\ref{defn:delta}) and (\ref{defn:gamma0}). From Proposition~\ref{prop:stable} we obtain two more constants $\ep_{1}$ and $M_{5}$. We also consider the constant $M_{1}$ of Proposition~\ref{prop:ultimate}, and the constant $r_{0}$ of Proposition~\ref{prop:asymptotic}. 

With a little abuse of notation, we consider the function 
\begin{equation}
E(u,v):=\frac{1}{2}|v|^{2}+\frac{1}{2}|Bu|^{2}-\frac{\lambda}{2}|A^{1/2}u|^{2}+\frac{1}{4}|A^{1/2}u|^{4},
\nonumber
\end{equation}
defined for every $(u,v)\in D(B)\times H$. In this way the classical energy $E(t)$ defined in (\ref{defn:E}) is just $E(u(t),u'(t))$.

The function $E(u,v)$ is continuous in $D(B)\times H$, and
\begin{equation}
E(\sigma_{0}e_{1},0)=-\frac{1}{4}\sigma_{0}^{4}{}<-\eta.
\nonumber
\end{equation}

As a consequence, there exists $\delta_{1}>0$ such that
\begin{equation}
|v|+|B(u-\sigma_{0}e_{1})|\leq\delta_{1}
\qquad\Longrightarrow\qquad
E(u,v)<-\eta.
\label{defn:delta-1}
\end{equation}

We can also assume that $\delta_{1}$ is small enough so that
\begin{equation}
|B(u-\sigma_{0}e_{1})|\leq\delta_{1}
\qquad\Longrightarrow\qquad
\langle u,Ae_{1}\rangle>0.
\label{defn:delta-1-bis}
\end{equation}

We claim that the conclusions of Theorem~\ref{thm:main} hold true if we choose
$$M_{0}:=\max\{M_{4},M_{5}\},$$
and we choose $\ep_{0}>0$ such that
\begin{equation}
\ep_{0}\leq\min\left\{1,\frac{\ep_{1}}{2},\frac{r_{0}}{2M_{0}},\frac{\delta_{1}}{2M_{0}}\right\}
\nonumber
\end{equation}
and
\begin{equation}
\ep_{0}^{2}\leq\frac{\gamma_{0}(1-\delta)\beta_{0}^{2}}{2M_{1}}.
\label{defn:ep0-2}
\end{equation}

\paragraph{\textmd{\textit{Alternative}}}

Let us assume that (\ref{hp:main}) is satisfied, and let $u(t)$ be any solution to (\ref{eqn:duffing-AB}). Let us set
\begin{equation}
L:=\limsup_{t\to +\infty}|u(t)|.
\nonumber
\end{equation}

We observe that $L$ is finite because of (\ref{th:ultimate-u}), and we distinguish two cases.

Let us assume that $L\leq\beta_{0}$. Since $\ep_{0}\leq 1$, we can apply Proposition~\ref{prop:unstable}, from which we deduce that in this case $u(t)$ satisfies (\ref{th:main-alternative}) with $\sigma=0$.
  
So it remain to consider the case $L>\beta_{0}$. In this case we claim that we are in the framework of Proposition~\ref{prop:stable} with $\eta$ given by (\ref{defn:beta}), namely there exists $T_{0}\geq 0$ for which the two inequalities in the upper box of (\ref{th:stable}) are satisfied.  

In order to check the first one, we observe that $\ep_{0}\leq\ep_{1}/2$, and therefore from assumption (\ref{hp:main}) it follows that $|f(t)|\leq\ep_{1}$ whenever $t$ is large enough.

In order to check the second one, we consider the function $\varphi(t):=\langle u(t),Pu(t)\rangle$. Due to (\ref{est:pu,u-u,u}), the function $\varphi(t)$ is bounded from above and
\begin{equation}
L_{1}:=\limsup_{t\to +\infty}\varphi(t)>(1-\delta)\beta_{0}^{2}.
\nonumber
\end{equation}

As a consequence, from Lemma~\ref{lemma:limsup} we deduce that there exists a sequence $t_{n}\to+\infty$ such that 
$$\varphi(t_{n})=\langle u(t_{n}),Pu(t_{n})\rangle\to L_{1}
\qquad\mbox{and}\qquad
\varphi'(t_{n})=\langle u'(t_{n}),Pu(t_{n})\rangle\to 0.$$

Now we observe that $E(t)=F(t)-\gamma_{0}\varphi(t)-2\gamma_{0}\varphi'(t)$. Setting $t=t_{n}$, and letting $n\to +\infty$, we obtain that
\begin{eqnarray*}
\limsup_{n\to +\infty}E(t_{n}) & = & \limsup_{n\to +\infty}F(t_{n})-\lim_{n\to +\infty}\varphi(t_{n})  \\[0.5ex]
& \leq & M_{1}\ep_{0}^{2}-\gamma_{0}(1-\delta)\beta_{0}^{2} \\[0.5ex]
& \leq & -2\eta,
\end{eqnarray*}
where in the inequalities we have exploited (\ref{th:ultimate}), the smallness assumption (\ref{defn:ep0-2}), and our definition (\ref{defn:beta}) of $\eta$. 

This shows that the two inequalities in the upper box of (\ref{th:stable}) are satisfied if we choose $T_{0}:=t_{n}$ with $n$ large enough. At this point, from the conclusion in the lower box of (\ref{th:stable}) we deduce that in this case $u(t)$ satisfies (\ref{th:main-alternative}) with $\sigma=\pm\sigma_{0}$.

\paragraph{\textmd{\textit{Asymptotic convergence}}}

Since $M_{0}\ep_{0}\leq r_{0}/2$, any pair of solutions satisfying~(\ref{th:main-alternative}) with the same $\sigma\in\{-\sigma_{0},0,\sigma_{0}\}$ satisfies also (\ref{hp:asymptotic}) with the same $\sigma$. At this point, (\ref{th:main-asymptotic}) follows from Proposition~\ref{prop:asymptotic}.

\paragraph{\textmd{\textit{Solutions in the stable regime}}}

Let us consider the case $\sigma=\sigma_{0}$ (but the argument is symmetric when $\sigma=-\sigma_{0}$). We claim that, when $f(t)$ satisfies (\ref{hp:main}) with our choice of $\ep_{0}$, the following characterization holds true: 
\begin{quote}
``a solution to (\ref{eqn:duffing-AB}) satisfies (\ref{th:main-alternative}) with $\sigma=\sigma_{0}$ if and only if there exists $T_{0}\geq 0$, possibly depending on the solution, for which the two inequalities in upper box of (\ref{th:stable}) hold true with $\eta$ given by (\ref{defn:beta}), and $\langle u(T_{0}),Ae_{1}\rangle>0$''.
\end{quote}

Let us prove this characterization. The ``if part'' is exactly Proposition~\ref{prop:stable}. As for the ``only if part'', it is enough to show that (\ref{hp:main}) and (\ref{th:main-alternative}) with $\sigma=\sigma_{0}$ imply that the two inequalities in the upper box of (\ref{th:stable}), and the further condition $\langle u(t),Ae_{1}\rangle>0$, hold true when $T_{0}$ is large enough. 

The first one follows from (\ref{hp:main}) because $\ep_{0}\leq\ep_{1}/2$. 

For the second one we observe that (\ref{th:main-alternative}) implies that
\begin{equation}
|u'(t)|+|B(u(t)-\sigma_{0}e_{1})|<2M_{0}\ep_{0}\leq\delta_{1}
\label{est:Bu-delta1}
\end{equation}
when $t$ is large enough, and hence from (\ref{defn:delta-1}) we deduce that $E(t)<-\eta$ for the same values of $t$.

The further condition holds true when $t$ is large enough because of (\ref{est:Bu-delta1}) and (\ref{defn:delta-1-bis}).

Given the characterization, we can prove our conclusions. Indeed, due to the continuous dependence on initial data, the set of initial data $(u_{0},u_{1})$ originating a solution $u(t)$ for which there required $T_{0}$ exists is an open set. In order to prove that it is nonempty, we choose $T_{0}\geq 0$ such that $|f(t)|\leq 2\ep_{0}\leq\ep_{1}$ for every $t\geq T_{0}$, and we consider the solution $u(t)$ to (\ref{eqn:duffing-AB}) with ``initial'' data
$$u(T_{0})=\sigma_{0}e_{1},
\qquad\quad
u'(T_{0})=0.$$

This solution fits in the assumptions of Proposition~\ref{prop:stable} because 
$$E(T_{0})=-\frac{1}{4}\sigma_{0}^{4}{}<-\eta.$$ 

We note that the extra condition $\langle u(T_{0}),Ae_{1}\rangle>0$ guarantees that the set of initial data for which the solution satisfies (\ref{th:main-alternative}) with $\sigma=\sigma_{0}$ is disjoint from the set of initial data for which the same relation is fulfilled with $\sigma=-\sigma_{0}$.

\paragraph{\textmd{\textit{Solutions in the unstable regime}}}

Due to the alternative of statement~(1), the set of initial data originating a solution satisfying (\ref{th:main-alternative}) with $\sigma=0$ is the complement of the set of initial data giving rise to solutions with $\sigma=\pm\sigma_{0}$. Since that set is the union of two open sets, the complement is necessarily closed, and nonempty because the phase space $\dath$ is connected and cannot be represented as the union of two disjoint nonempty open sets.\qed 

\setcounter{equation}{0}
\section{The concrete case (proof of Proposition~\ref{prop:concrete})}\label{sec:concrete}

We need to check that the operators $A$ and $B$ satisfy all the requirements in Definition~\ref{defn:gap}. It is a classical result that $A$ is a self-adjoint operator, and it satisfies (\ref{defn:mu-2}) with $\mu_{2}:=\pi^{2}$. Indeed, in this concrete case it turns out that $D(A^{1/2})=H^{1}_{0}((0,1))$, and (\ref{defn:mu-2}) reduces to
$$\int_{0}^{1}u_{x}(x)^{2}\,dx\geq\pi^{2}\int_{0}^{1}u(x)^{2}\,dx
\qquad
\forall u\in H^{1}_{0}((0,1)),$$
which is Poincaré inequality.

The operator $B$ is the square root of $B^{2}$, and its domain is
$$D(B)=H^{2}_{0}((0,1))\subseteq H^{2}((0,1))\cap H^{1}_{0}((0,1))=D(A),$$
as required. Moreover, inequality (\ref{defn:mu-1}) holds true with $\mu_{1}:=1$ because
$$|Bu|^{2}=\langle B^{2}u,u\rangle=|Au|^{2}=\int_{0}^{1}u_{xx}(x)^{2}\,dx
\qquad
\forall u\in H^{2}_{0}((0,1)).$$

So it remains to check that (\ref{hp:e1}) and (\ref{hp:e1-perp}) hold true with the values of $\lambda_{1}$ and $\lambda_{2}$ given in the statement, and with a suitable $e_{1}$. To this end, we begin by investigating all nontrivial solutions to (\ref{hp:e1}), and then we conclude the proof in two alternative ways.

\paragraph{\emph{\textmd{Nontrivial solutions to equation (\ref{hp:e1})}}}

We look for all pairs $(\lambda,\varphi)$, where $\lambda$ is a positive real number and $\varphi$ is a smooth function that satisfies
\begin{equation}
\varphi_{xxxx}(x)+\lambda\varphi_{xx}(x)=0
\qquad
\forall x\in(0,1),
\label{ode}
\end{equation}
with boundary conditions 
\begin{equation}
\varphi(0)=\varphi_{x}(0)=\varphi(1)=\varphi_{x}(1)=0.
\label{bc}
\end{equation}

Setting $\sqrt{\lambda}=:2\alpha$ for the sake of shortness, all solution to (\ref{ode}) are of the form
$$ \varphi(x)= c_1 + c_2\,x + c_3 \cos\left(2\alpha x\right) + c_4 \sin\left(2\alpha x\right).$$

From the boundary conditions in $x=0$ we deduce that $c_1+c_3=c_2+2\alpha c_4= 0$, so that we can restrict to solutions of the form
$$ \varphi(x)= c_{1}\left[ 1- \cos\left(2\alpha x\right) \right] + c_{4}\left[ \sin\left(2\alpha x\right)- 2\alpha x \right].$$

The boundary conditions in $x=1$ are now equivalent to the system
\begin{equation}
\left\{\begin{array}{l}
[1-\cos(2\alpha)]\cdot c_{1}+[\sin(2\alpha)- 2\alpha]\cdot c_{4}= 0, \\[1ex]
2\alpha\sin(2\alpha)\cdot c_{1}+2\alpha[\cos(2\alpha)-1]\cdot c_{4}=0.
\end{array}\right.
\nonumber
\end{equation}

With some standard algebra, we can show that this system has a nontrivial solution if and only if
$$\sin^{2}\alpha=\alpha\sin\alpha\cos\alpha.$$ 

The solutions to this equation are all solutions to $\sin\alpha=0$, namely the values of the form $\alpha=k\pi$ (with $k$ any positive integer), and all solutions to $\tan\alpha=\alpha$, which are an infinite sequence $\alpha_{k}$, with one element in each interval of the form $(k\pi,(k+1/2)\pi)$. Recalling that $\lambda=4\alpha^{2}$, the required eigenvalues are those of the form $4\pi^{2}k^{2}$, with corresponding eigenfunctions
$$\varphi_{k}(x)=1-\cos(2k\pi x),$$
and those of the form $4\alpha_{k}^{2}$, with corresponding eigenfunctions of the form
$$\varphi_{k}(x)=\alpha_{k}(1-\cos(2\alpha_{k}x))+\sin(2\alpha_{k}x)-2\alpha_{k}x.$$

\paragraph{\textmd{\emph{Conclusion through variational approach}}}

We show that (\ref{hp:e1}) and (\ref{hp:e1-perp}) hold true with $\lambda_{1}=4\pi^{2}$, and consequently $e_{1}$ equal to a suitable multiple of $1-\cos(2\pi x)$, and $\lambda_{2}=4\alpha_{1}^{2}$. To this end, we consider the minimum problem (\ref{defn:lambda-1}), which now reads as
$$\min\left\{\int_{0}^{1}\varphi_{xx}(x)^{2}\,dx:\varphi\in H^{2}_{0}((0,1))\mbox{ and } \int_{0}^{1}\varphi_{x}(x)^{2}\,dx=1\right\},$$
and the minimum problem (\ref{defn:lambda-2}), which now reads as
$$\min\left\{\int_{0}^{1}\varphi_{xx}(x)^{2}\,dx:\varphi\in H^{2}_{0}((0,1)),\ \int_{0}^{1}\varphi_{x}(x)^{2}\,dx=1,\ \int_{0}^{1}\sin(2\pi x)\varphi_{x}(x)\,dx=0\right\}.$$

In both cases a standard application of the direct method in the calculus of variations shows that the minimum exists, and any minimizer satisfies (\ref{ode}) with boundary conditions (\ref{bc}), and $\lambda$ equal to the minimum value. It follows that the two minimum values are the two smallest values of $\lambda$ for which (\ref{ode})--(\ref{bc}) has nontrivial solutions, and from the previous analysis we know that these values are exactly $4\pi^{2}$ and $4\alpha_{1}^{2}$. This proves (\ref{hp:e1}) and (\ref{hp:e1-perp}) as required.

\paragraph{\textmd{\emph{Conclusion through functional analytic approach}}}

As discussed in section~\ref{appendix:concr} of the appendix, the operator $A^{-1}B^{2}$ is symmetric with compact inverse in $H^{1}_{0}((0,1))$, and therefore the eigenfunctions that we found above span a dense subspace of $L^{2}((0,1))$, and they are orthogonal with respect to the scalar product of $H^{1}_{0}((0,1))$. At this point (\ref{hp:e1}) and (\ref{hp:e1-perp}) follow from the classical variational characterization of the two smallest eigenvalues, applied in this case to the operator $A^{-1}B^{2}$ in the space $H^{1}_{0}((0,1))$.


\appendix

\setcounter{equation}{0}
\section{Appendix}

In the second paragraph of section~\ref{sec:statements} we described in finite dimension the role of the eigenvalues of $A^{-1}B^{2}$ in the study of the negative inertia index of $B^{2}-\lambda A$ as a function of $\lambda$. In the rest of the paper we developed our theory in the infinite dimensional case without mentioning $A^{-1}B^{2}$ explicitly. 

In this appendix we present a possible functional setting in which the spectral theory can be applied to the operator $A^{-1}B^{2}$, both in the general and in the concrete case.

\subsection{The correct framework for $A^{-1}B^2$ in general} 

Let $H$ be a Hilbert space, and let $A$ and $B$ be two coercive self-adjoint unbounded operators on $H$ with dense domains $D(B)\subseteq D(A)$. Then a reasonable definition of $A^{-1}B^{2}$ seems to be the following.

Let us consider the Hilbert spaces $V:=D(A^{1/2})$ and $W:=D(B)$. If we identify $H$ with its dual space $H'$, then we have the inclusions
$$W\subseteq V\subseteq H=H'\subseteq V'\subseteq W'.$$

With these notations we can consider $A$ as a bounded operator $\widehat{A}:V\to V'$, and represent the scalar product in $V$ in terms of the duality pairing as 
$$\langle u,v\rangle_{V}=\langle \widehat{A}u,v\rangle_{V',V}
\qquad
\forall(u,v)\in V^{2}.$$

Similarly, we con consider $B$ as a bounded operator $\widehat{B}:W\to H$, whose adjoint is a bounded operator $\widehat{B}^{*}:H'\to W'$ with $\widehat{B}^{*}u=Bu$ if $u\in W\subseteq H'$.

Now we can consider the unbounded operator $C$ in $V$ with domain
\begin{equation}
D(C):=\left\{u\in W:\widehat{B}^{*}Bu\in V'\right\},
\nonumber
\end{equation}
and defined by
\begin{equation}
Cu:=\widehat{A}^{-1}\widehat{B}^{*}Bu
\qquad
\forall u\in D(C).
\nonumber
\end{equation}

We claim that $C$ is an extension of $A^{-1}B^{2}$ that is symmetric and maximal monotone as an unbounded operator in $V$.

To begin with, for every $u$ and $v$ in $D(C)$ it turns out that
$$\langle v,Cu\rangle_{V}=\langle v,\widehat{A}Cu\rangle_{V',V}=\langle v,\widehat{B}^{*}Bu\rangle_{V',V}=\langle \widehat{B}v,Bu\rangle_{H}=\langle Bv,Bu\rangle_{H},$$
which proves that $C$ is symmetric and monotone.

It remain to show that $C$ is maximal, namely that, for every $f\in V$, the equation $u+Cu=f$ has a (unique) solution $u\in D(C)$. Applying $\widehat{A}$ to both sides, this equation becomes
$$\widehat{A}u+\widehat{B}^{*}Bu=\widehat{A}f\in V'.$$

Now the operator $\widehat{A}+\widehat{B}^{*}B$ is coercive from $W$ to $W'$, and hence surjective (see for example~\cite[Corollary~14]{brezis68}). Since $\widehat{A}f\in V'$, the solution belongs to $D(C)$.


\subsection{The operator $A^{-1}B^2$ in the concrete case}\label{appendix:concr}

Instead of fussing with generalities, we give an explicit description of the operator $A^{-1}B^{2}$ in the case where $A$ and $B$ are as in Proposition~\ref{prop:concrete}.

Let us consider the Hilbert space $V=H^{1}_{0}((0,1))$ with scalar product
\begin{equation}
\langle u,v\rangle_{V}:=\int_{0}^{1}u_{x}(x)v_{x}(x)\,dx
\qquad
\forall (u,v)\in H^{1}_{0}((0,1))^{2}.
\label{defn:scal-prod-V}
\end{equation}

Let us consider the unbounded linear operator $C$ in $V$ with domain
\begin{equation}
D(C):=H^{3}((0,1))\cap H^{2}_{0}((0,1)),
\nonumber
\end{equation}
and defined by
\begin{equation}
[Cu](x):=-u_{xx}(x)+u_{xx}(0)+[u_{xx}(1)-u_{xx}(0)]x
\qquad
\forall u\in D(C).
\label{defn:Cu}
\end{equation}

We observe that
\begin{equation}
Cu=A^{-1}B^{2}u
\qquad
\forall u\in D(B^{2})=H^{4}((0,1))\cap H^{2}_{0}((0,1)),
\nonumber
\end{equation}
and therefore $C$ is a natural extension of $A^{-1}B^{2}$. Indeed $B^{2}u=u_{xxxx}\in L^{2}((0,1))$, and hence $A^{-1}B^{2}u$ is the solution $v(x)$ to equation $-v_{xx}=u_{xxxx}$ in $(0,1)$, with Dirichlet boundary conditions in $x=0$ and $x=1$. The solution is exactly the function $Cu$ defined in (\ref{defn:Cu}).

We claim that $C$ is a symmetric positive operator in $V$ with domain $D(C)$, and compact inverse. If we prove this claim, then the eigenfunctions of $C$ are a basis of $H^{1}_{0}((0,1))$, and hence also a basis of $L^{2}((0,1))$, but orthogonal with respect to the scalar product (\ref{defn:scal-prod-V}). These eigenfunctions are exactly the solutions to (\ref{ode}) that we characterized in section~\ref{sec:concrete}.

To begin with, for every $u$ and $v$ in $D(C)$ it turns out that
\begin{eqnarray*}
\langle v,Cu\rangle_{V} & = & \int_{0}^{1}v_{x}(x)\cdot[Cu]_{x}(x)\,dx   \\
& = & \int_{0}^{1}v_{x}(x)\cdot\{-u_{xxx}(x)+u_{xx}(1)-u_{xx}(0)\}\,dx  \\
& = & -\int_{0}^{1}v_{x}(x)\cdot u_{xxx}(x)\,dx  \\
& = & \int_{0}^{1}v_{xx}(x)\cdot u_{xx}(x)\,dx,
\end{eqnarray*}
which is enough to conclude that $C$ is both symmetric and positive.

It remains to show that the inverse is compact, namely that for every $f\in H^{1}_{0}((0,1))$ the equation
\begin{equation}
-u_{xx}(x)+u_{xx}(0)+[u_{xx}(1)-u_{xx}(0)]x=f(x)
\qquad
\forall x\in(0,1)
\label{defn:C-1}
\end{equation}
has a unique solution $Tf\in D(C)$ (note that lying in this domain entails four boundary conditions), and $T$ is compact as an operator $T:H^{1}_{0}((0,1))\to H^{1}_{0}((0,1))$. 

Uniqueness follows from the positivity of $C$. Existence follows from the explicit formula for the solution. Indeed, if we set
$$F(x):=\int_{0}^{x}f(t)\,dt,
\qquad\qquad
\widehat{F}(x):=\int_{0}^{x}F(t)\,dt,$$
then a standard computation shows that the solution to (\ref{defn:C-1}) is
\begin{equation}
u(x):=-\widehat{F}(x)+\left(3\widehat{F}(1)-F(1)\right)x^{2}+\left(F(1)-2\widehat{F}(1)\right)x^{3}.
\nonumber
\end{equation}

The same formula reveals that if a sequence $\{f_{n}\}$ is bounded in $H^{1}_{0}((0,1))$, then the sequence of corresponding solutions is bounded in $H^{3}((0,1))$, and therefore relatively compact in $H^{1}_{0}((0,1))$.


\subsubsection*{\centering Acknowledgments}

The first two authors are members of the ``Gruppo Nazionale per l'Analisi Matematica, la Probabilit\`{a} e le loro Applicazioni'' (GNAMPA) of the ``Istituto Nazionale di Alta Matematica'' (INdAM). 


\label{NumeroPagine}

\end{document}